


%
%
%


\input amstex

\documentstyle{amsppt}
\magnification=\magstep1
\NoRunningHeads
\NoBlackBoxes
\TagsOnRight

\input epsf
\def\seteps#1#2#3#4{\vskip0in\relax\noindent\hskip#1\relax
 \epsfxsize=#2\epsfysize=#3\epsfbox{#4}}
\def\centereps#1#2#3{\vskip0in\relax\centerline{\epsfxsize=#1\epsfysize=#2\epsfbox{#3}}}

\define\qqand{\qquad\text{and}\qquad}
\define\qand{\quad\text{and}\quad}
\define\eg{{\it e.g.}}
\define\ie{{\it i.e.}}
\define\cf{{\it cf.\ }}

\define \({\left(}  
\define \){\right)_{\vphantom{\leavevmode\lower1pt\hbox{x}}}}
\define \der#1#2{\frac{d#1}{d#2}}        
\define\Perp{{}_{\!\sssize\perp}}        
\define \<#1>{\left<#1\right>}           

\define\del{\nabla\futurelet\next\negthinspaceifsubscript}  
\define\negthinspaceifsubscript{\ifx\next_\!\fi}

\define\grad{\operatorname{grad}}
\redefine\l{\ell}
\define\a{\alpha}
\redefine\b{\beta}
\define\g{\gamma}
\redefine\phi{\varphi}
\define\phit{\tilde\phi}
\define\bdy{\partial}

\redefine\P{\Cal P}
\define\kg{\kappa_g}
\define\E#1{\Bbb E^{#1}}
\redefine\H#1{\Bbb H^{#1}}
\let\section=\S  
\redefine\S#1{\Bbb S^{#1}}
\define\bK{{\bdy K}}
\define\ambientK{\bold K}

\define\Line#1{\overleftrightarrow{#1}}
\define\ray#1{\overrightarrow{#1}}
\define\seg#1{\overline{#1}}
\redefine\D#1{\Bbb D^{#1}}

\document

\topmatter
\title Integral-Geometric Formulas for Perimeter in $\S2$, $\H2$, and Hilbert Planes\endtitle
\author Ralph Alexander, I.\,D.\,Berg\\
{\rm University of Illinois at Urbana--Champaign}\\
\smallskip
Robert L. Foote \\
{\rm Wabash College}
\endauthor
\address Department of Mathematics, Urbana, IL 61801 \endaddress
\email jralex\@math.uiuc.edu,\quad berg\@math.uiuc.edu \endemail
\address Department of Mathematics \& Computer Science, Crawfordsville, IN 47933 \endaddress
\email footer\@wabash.edu \endemail
\abstract We develop two types of integral formulas for the perimeter of a convex body $K$ in planar geometries. We derive Cauchy-type formulas for perimeter in planar Hilbert geometries. Specializing to $\H2$ we get a formula that appears to be new. In the projective model of $\H2$ we have $\P = \frac12 \int w\,d\phi$. Here $w$ is the Euclidean length of the projection of $K$ from the ideal boundary point $R=(\cos\phi,\sin\phi)$ onto the diametric line perpendicular to the radial line to $R$ (the image of $K$ may contain points outside the model). We show that the standard Cauchy formula $\P = \int \sinh r\,d\omega$ in $\H2$ follows, where $\omega$ is a central angle perpendicular to a support line and $r$ is the distance to the support line.

The Minkowski formula $\P = \int \kg \rho^2 \,d\theta$ in $\E2$ generalizes to
$\P = \frac1{4\pi^2} \int \kg L(\rho)^2\, d\theta + \frac k{2\pi} \int A(\rho)\,ds$ in $\H2$ and $\S2$. Here $(\rho,\theta)$ and $\kg$ are, respectively, the polar coordinates and geodesic curvature of $\bK$, $k$ is the (constant) curvature of the plane, and $L(\rho)$ and $A(\rho)$ are, respectively, the perimeter and area of the disk of radius $\rho$. In $\E2$ this is locally equivalent to the Cauchy formula $\P = \int r\,d\omega$ in the sense that the integrands are pointwise equal. In contrast, the corresponding Minkowski and Cauchy formulas are not locally equivalent in $\H2$ and $\S2$.
\endabstract
\dedicatory Dedicated to the memory of our friend and colleague Felix Albrecht. \enddedicatory
\thanks The authors would like to express thanks to D. Chakerian, R. Schneider, J. Sullivan, P. Tondeur, and the referee for their helpful suggestions.
\endthanks
\thanks To appear in the Rocky Mountain Journal of Mathematics \endthanks
\endtopmatter

\document
\baselineskip15pt




\heading 1. Introduction\endheading

\subheading{Overview}
There are at least two natural integral-geometric approaches relating the perimeter of a convex body $K$ in $\E2$ to its other geometric properties. There is the beautiful Cauchy formula 
$$
\dsize \P = \frac{1}{2}\int_0^{2\pi}w\,d\phi, \tag{1.1}
$$
where $w$ is the width of the body at the points where the support lines are at angle $\phi$ to some fixed axis (see Figure~1). Two variations of this are the one-sided formulas
$$
\P = \int_0^{2\pi} r\,d\phi \qqand \P = \int_0^{2\pi} r\,d\omega, 
\tag{{1.2},{1.3}}
$$
where $r$ is the ``support'' function, that is, the signed distance from the origin to the ``right'' support line as viewed from infinity in the direction given by $\phi$, and $\omega$ is the central angle from the fixed axis to the ray from the origin to the closest point on the right support line (here $\omega = \phi + \pi/2$). There is also the Minkowski approach, which has two forms,
$$
\dsize \P = \int_\bK r\kg \,ds = \int_\bK \rho^2 \kg\,d\theta, \tag{1.4}
$$
where $\rho$ and $\theta$ are the polar coordinates of a point on the boundary curve $\bK$ and $\kg$ is the geodesic curvature of $\bK$. We call these the ``support'' and ``polar'' forms, respectively. (Note that $\frac12r\,ds$ and $\frac12\rho^2\,d\theta$ are equal since both are equal to $dA$, the area of an infinitesimal triangle with one vertex at the origin and two vertices infinitesimally close on $\bK$, or from considering similar triangles as in Figure~7.)

\smallskip
\centereps{2.5in}{2in}{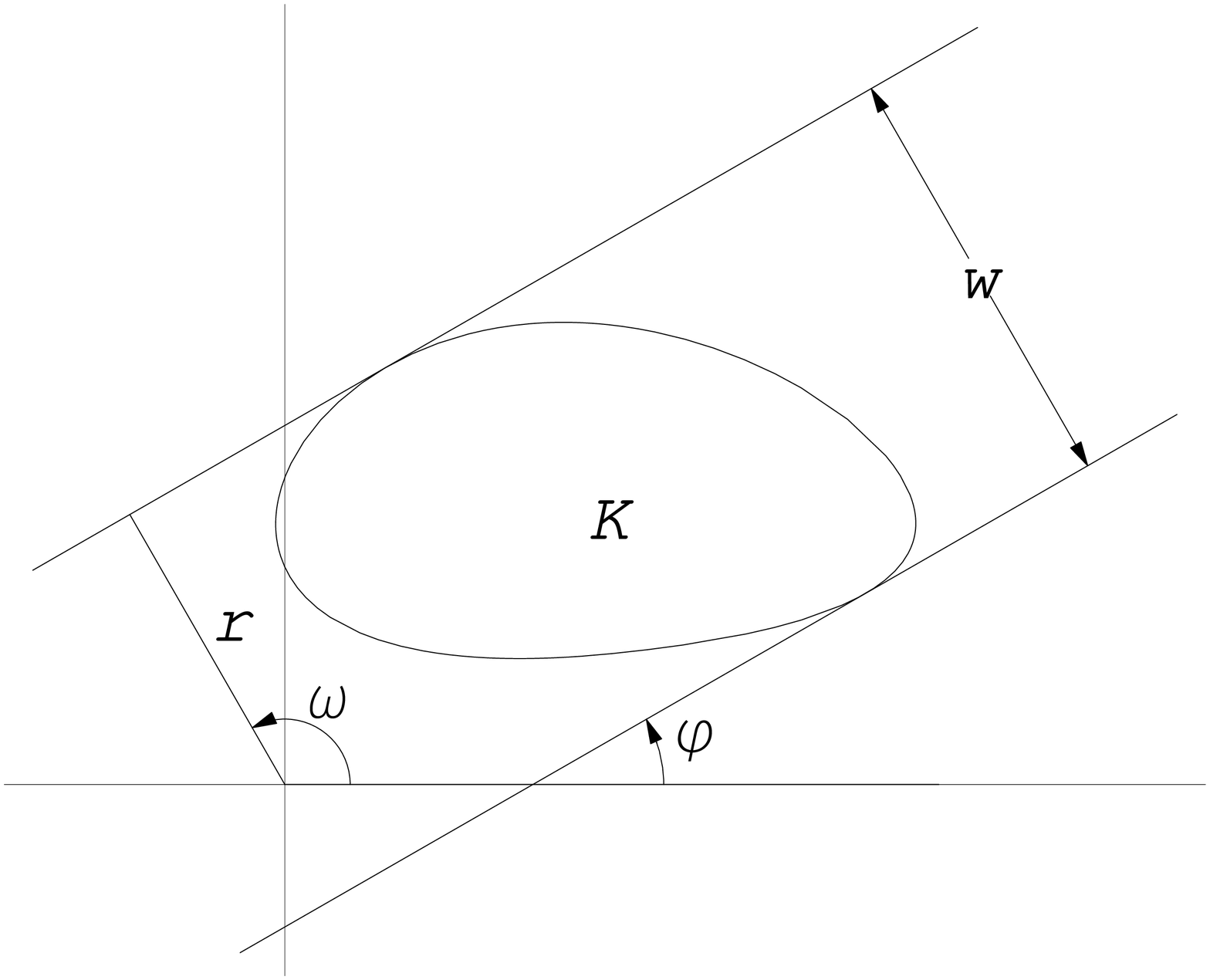}

\centerline{Figure~1.}
\smallskip


Both of these approaches have natural generalizations to $\H2$ and $\S2$. In \section 2 we generalize the polar Minkowski formula. In \section 3 we generalize the Cauchy formulas to two-dimensional Hilbert geometries. In \section4 we specialize our Cauchy formulas to $\H2$. In $\E2$ the Minkowski and Cauchy formulas are locally equivalent in the sense that the integrands in \thetag{1.3} and \thetag{1.4} are equal at each point of $\bK$ (since $\kg\,ds = d\omega$). Somewhat surprisingly this equivalence fails in $\H2$ and $\S2$, as we show in \section5 where we reformulate the Cauchy form, allowing comparison with the Minkowski form. In \section 6 we develop formulas relating the various measures that arise on $\bK$.

The Cauchy formulas and our generalizations of them \thetag{1.6, 1.7, 1.8, 1.9} hold when $K$ is convex. In contrast, the polar Minkowski formula and our generalization of it \thetag{1.5} hold when $\bK$ is $C^2$, with no assumption of convexity. In fact, they compute the length of any closed, $C^2$ curve.

Our general references for integral geometry are \cite{S5,So}.

\subheading{Minkowski Formulas}
The polar form of the Minkowski formula \thetag{1.4} was posed as a problem by Antonio Montes in the {\sl American Mathematical Monthly}, where it attracted several interesting solutions \cite{M}. It generalizes easily to $\E n$. The support form in $\E3$, as described in Burago and Zalgaller \cite{B\&Z, pg.~51}, is (with slight modification) $A = \int_\bK rh\,dA$, where $K$ is a compact, convex set with smooth boundary, $A$ is the area of $\bK$, $h$ is the mean curvature of $\bK$ at a point $p$, and $r$ is the signed distance from the origin to the support plane at $p$. (See also \cite{S6, pg.~326}.) Similar to the case in $\E2$, one can use $r\,dA = \rho^3\,d\Theta = 3\,dV$ to convert this to the polar form $A = \int_\bK \rho^3 h\,d\Theta$, where $\rho$ is the distance from the origin to $p$ and $d\Theta$ is the solid angle element with respect to the origin. None of the solvers of the problem posed by Montes remarked on the Minkowski result in $\E3$ or happened to mention Minkowski's name. It seems appropriate by reason of priority to refer to this result by Minkowski's name although each solver offered his own insights. The second author (Berg) became interested in this area by considering Montes' problem and acknowledges valuable correspondence with Montes.

The polar Minkowski formula \thetag{1.4} also generalize nicely to spaces of constant curvature $k$. In \section 2 we show that
$$
\P = \int_\bK \left(\frac{L(\rho)}{2\pi}\right)^2\! \kg \,d\theta 
    + k\int_\bK \frac{A(\rho)}{2\pi}\,ds  \tag{1.5}
$$
in $\H2$ and $\S2$, where $\rho$ and $\theta$ are intrinsic polar coordinates with respect to some fixed point, $L(\rho)$ is the circumference of a circle of radius $\rho$, and $A(\rho)$ is the area of a disk of radius $\rho$ (formulas for these are in the appendix). The proof of \thetag{1.5} shows that the integrand can be rewritten as $ds + d\(\frac{L(\rho)}{2\pi} \sin\a\)$, where $\a$ is the angle between the radial and normal vectors (see Figures~2 and~3). 

\subheading{Cauchy Formulas}
The natural generalizations of the Cauchy formulas \thetag{1.1, 1.2, 1.3} to the hyperbolic plane pass through results for Hilbert geometries. They take particularly satisfying forms in the Beltrami-Klein model, which we establish in \section4. To discuss these, we introduce several related quantities, illustrated in Figure~2.

The disk represents the Beltrami-Klein model of $\H2$. Let $K$ be any convex body in $\H2$. Let $R$ be an arbitrary point at infinity, determined by its central angle $\phi$. The line $\Line{OS}$ through the origin $O$ perpendicular to the terminal ray $\ray{OR}$ of $\phi$ will be called the $\phi$-{\sl normal}. Consider the two support lines of $K$ passing through $R$. In particular, note where they intersect the $\phi$-normal, extending, if necessary, past the boundary of the model. Let $w$ denote the Euclidean distance between these intersections, that is, $w$ is the Euclidean length of the projection of $K$ from $R$ onto the $\phi$-normal. 

\centereps{3in}{3in}{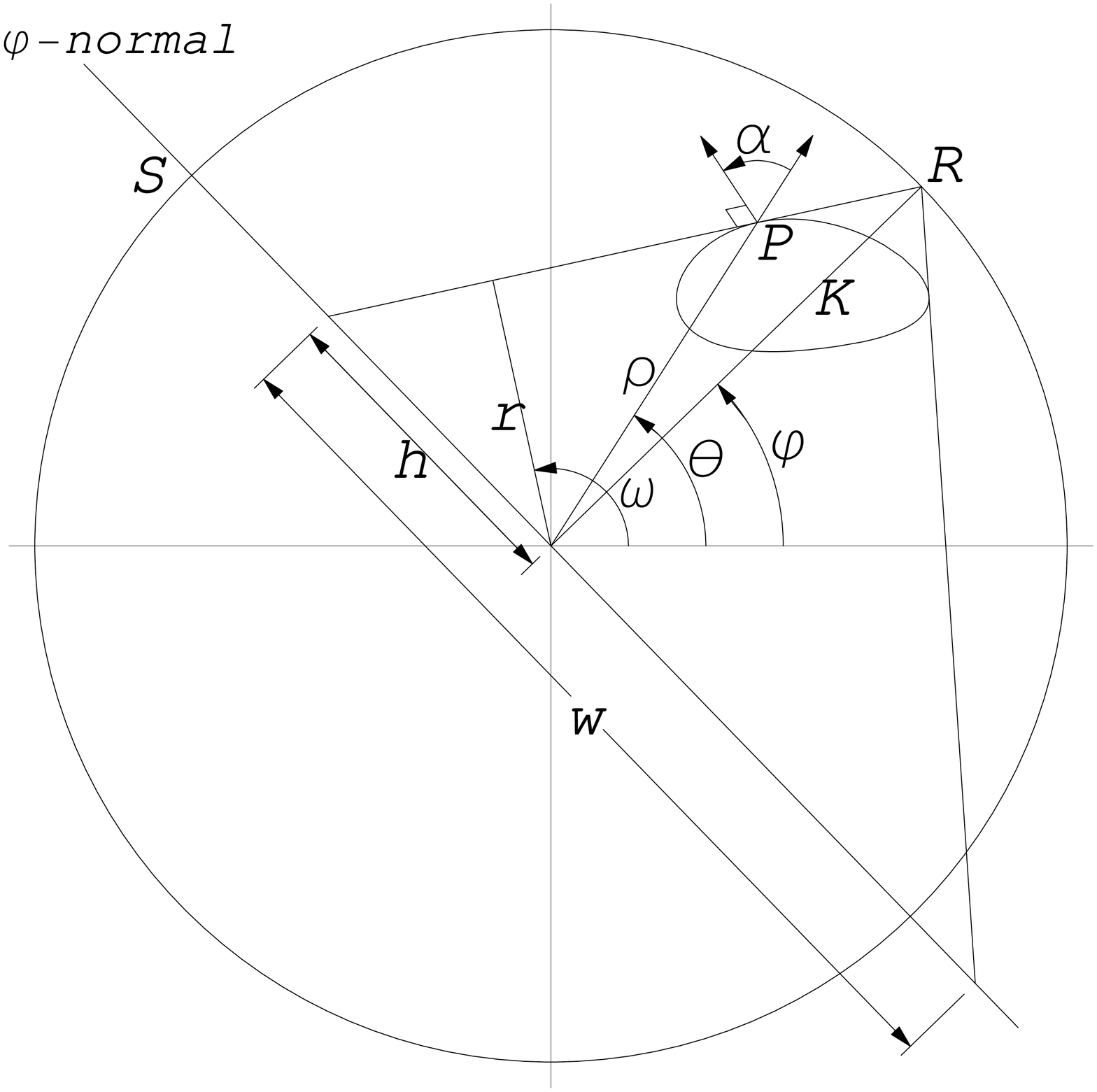}

\centerline{Figure~2. Beltrami-Klein model of $\H2$.}
\centerline{Lengths $\rho$ and $r$ are hyperbolic; $w$ and $h$ are Euclidean.}
\bigskip

Our generalization of \thetag{1.1} in this setting (see \section3) is
$$
\P = \frac{1}{2} \int_0^{2\pi} w\,d\phi. \tag{1.6}
$$
In this sense $w$ is the ``width'' of the body viewed from $\phi$. From this we derive a generalization of \thetag{1.2}, namely the one-sided formula
$$
\dsize \P = \int_0^{2\pi} h\,d\phi, \tag{1.7}
$$
where $h$ is the signed Euclidean length of the projection of the right hand support line as viewed from $R$. We call \thetag{1.6} and \thetag{1.7} our {\sl projective Cauchy formulas.} They provide an appealing characterization of length in $\H2$ and appear to be unpublished. Santal\'o's thorough treatment of Crofton's formulas (\eg\ \cite{S5}) has unplumbed depths but \thetag{1.6} and \thetag{1.7} do not seem to be on the surface, and their Hilbert geometry analogue \thetag{1.10} less so.

From the projective Cauchy formula \thetag{1.7} we deduce generalizations of \thetag{1.2} and \thetag{1.3},
$$
\dsize \P=\int_0^{2\pi} \frac{L(r)}{2\pi}d\phi \qqand 
\P=\int_0^{2\pi} \frac{L(r)}{2\pi}d\omega, \tag{1.8, 1.9}
$$
where (see Figure~2) $r$ is the signed {\sl hyperbolic\/} distance to the right hand support line, $\omega$ is the central angle of that line's closest point to the origin, and, as in \thetag{1.5}, $L(r)$ is the hyperbolic circumference of the circle of radius $r$ (when $r$ is negative, $L(r)$ is also taken to be negative). These forms are intrinsic. Furthermore, \thetag{1.9} applies to $\H2$, $\E2$ and $\S2$, and is known, \cf Santal\'o \cite{S2, S3, S4}. (In $\S2$, $K$ must lie in the open hemisphere centered at the origin.)  The integrands of \thetag{1.6} and \thetag{1.7} can also be viewed as intrinsic, although they may appear to depend on the Beltrami-Klein model: $h(\phi)$ is the signed inversive product of the right support line and the $\phi$-normal. Similarly, $w(\phi)$ is the sum of two inversive products. Note that while \thetag{1.1}, \thetag{1.2}, and \thetag{1.3} are obviously equivalent, the relationship between \thetag{1.6}, \thetag{1.7}, \thetag{1.8}, and \thetag{1.9} is less evident.

Finishing the quantities indicated in Figure~2, $P$ is a point of contact of $K$ with the right support line, $\rho$ and $\theta$ are the intrinsic polar coordinates of $P$ (used in \thetag{1.5}), and $\a$ is the intrinsic angle from the radial vector to the hyperbolic outward normal to $\bK$ at $P$. Note that if every support line of $K$ meets $\bK$ in a unique point, that is, $\bK$ contains no geodesic segments, then all of the quantities are uniquely determined by $\phi$ and $\omega$. If every point $P$ in $\bK$ has a unique support line, that is, $\bK$ has no corners, then all of the quantities are uniquely determined by $P$. Finally, if the origin is contained in the interior of $K$, then $P$ is uniquely determined by $\theta$. Figure~3 illustrates the portions of Figure~2 that are valid in all three geometries.

\centereps{3in}{2in}{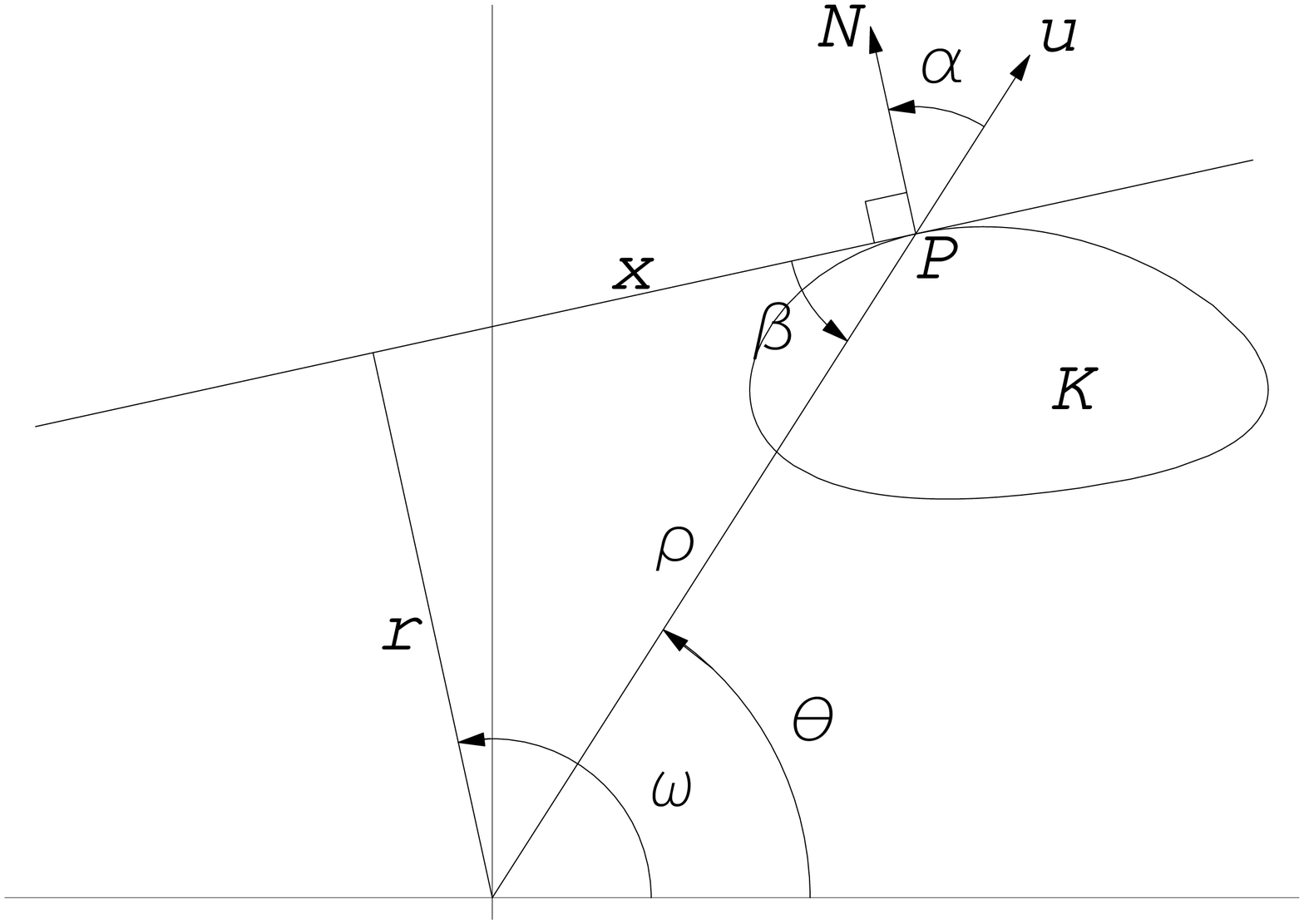}
\centerline{Figure~3.}

Our formula \thetag{1.6} in the hyperbolic plane is a special case of a Cauchy formula in two-dimensional Hilbert geometry, which we derive in \section3. For any bounded, convex body $\ambientK$ in the plane, Hilbert \cite{H, Appendix 1} showed that we get a metric defining a Hilbert distance $h$ on the interior of $\ambientK$, by letting 
$h(P, Q) = \frac12\log \(\frac{AQ}{AP}\,\frac{BP}{BQ}\)$, as in Figure~4. The unit disk yields $\H2$ with curvature $k=-1$. In this geometry the geodesics are Euclidean line segments. Let $K$ be a compact, convex body in the interior of $\ambientK$. The Cauchy formula for the Hilbert perimeter of $K$ is
$$
\P = \frac12 \int_0^{2\pi}
        \bigl(\cot\psi_1(\phi) - \cot\psi_2(\phi)\bigr)\,d\phi, \tag{1.10}
$$
which becomes \thetag{1.6} in $\H2$. Here (see Figure~5) $\phi$ is taken to be the angle from a fixed axis to the outward normal of a support line of $\ambientK$ at some point $R$. The angles $\psi_1$ and $\psi_2$ are taken from the support line to, respectively, the right and left support lines of $K$ as viewed from and passing though $R$ (so that $0 < \psi_1 \le \psi_2 < \pi$). This integral makes sense since $\psi_1$ and $\psi_2$ are defined for almost every $\phi$ in $[0,2\pi]$.  Formula \thetag{1.10} follows from the special case that gives the Hilbert distance between points $P$ and $Q$ in $\ambientK$, namely 
$2h(P,Q) = 
\frac12 \int_0^{2\pi}\bigl(\cot\psi_1(\phi) - \cot\psi_2(\phi)\bigr)\,d\phi$.
Literally this is the perimeter of the two-gon with vertices $P$ and $Q$ in Figure~6. We establish this formula for $h$ by direct calculation; we discovered it, however, by pondering a calculation of cross-ratios in polygons in the first author's (R. Alexander) work on Hilbert geometries \cite{A}.


\vskip.2in
\seteps{-.25in}{2in}{2in}{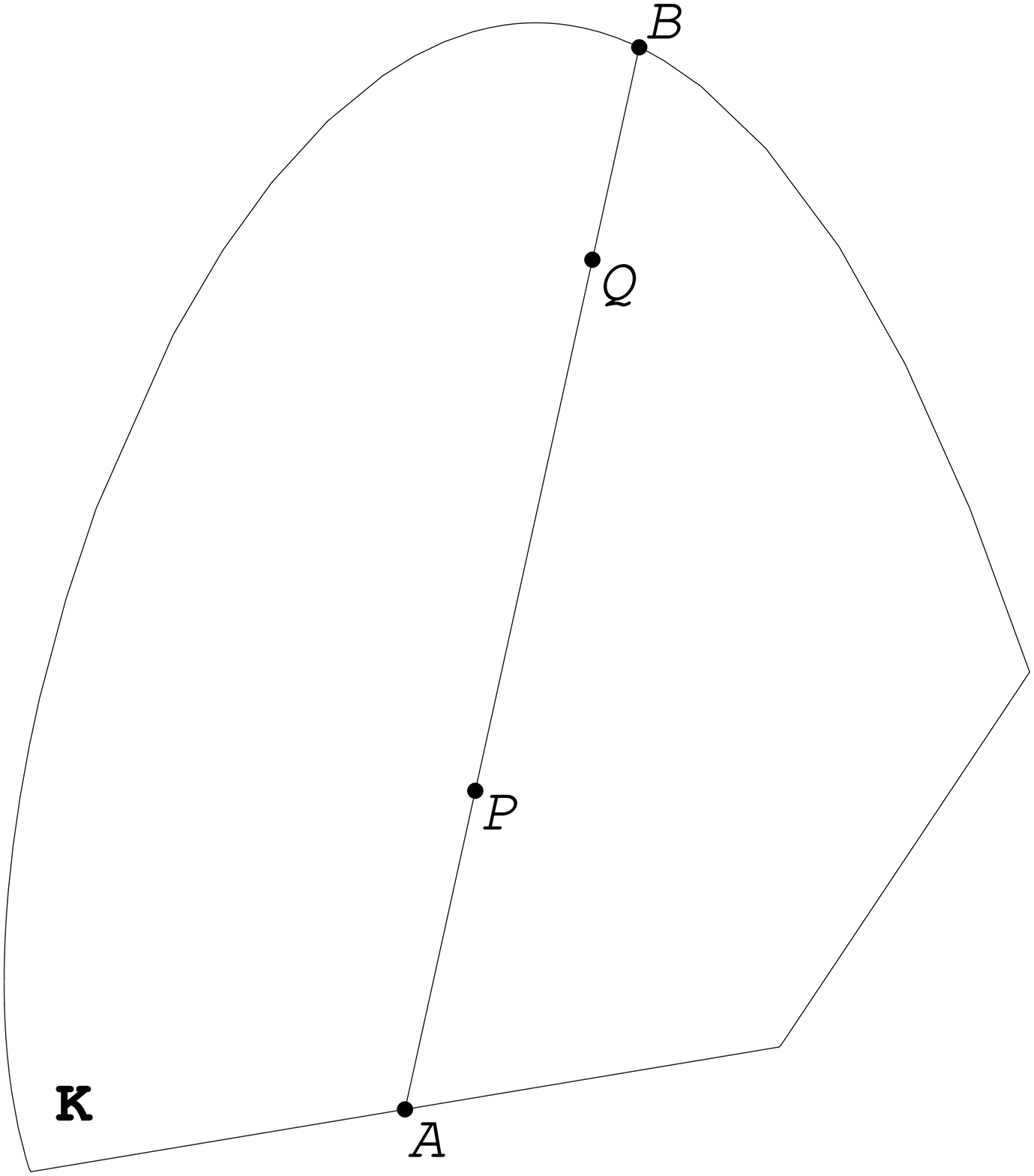}
\vskip-2.2in
\seteps{1.75in}{2in}{2.2in}{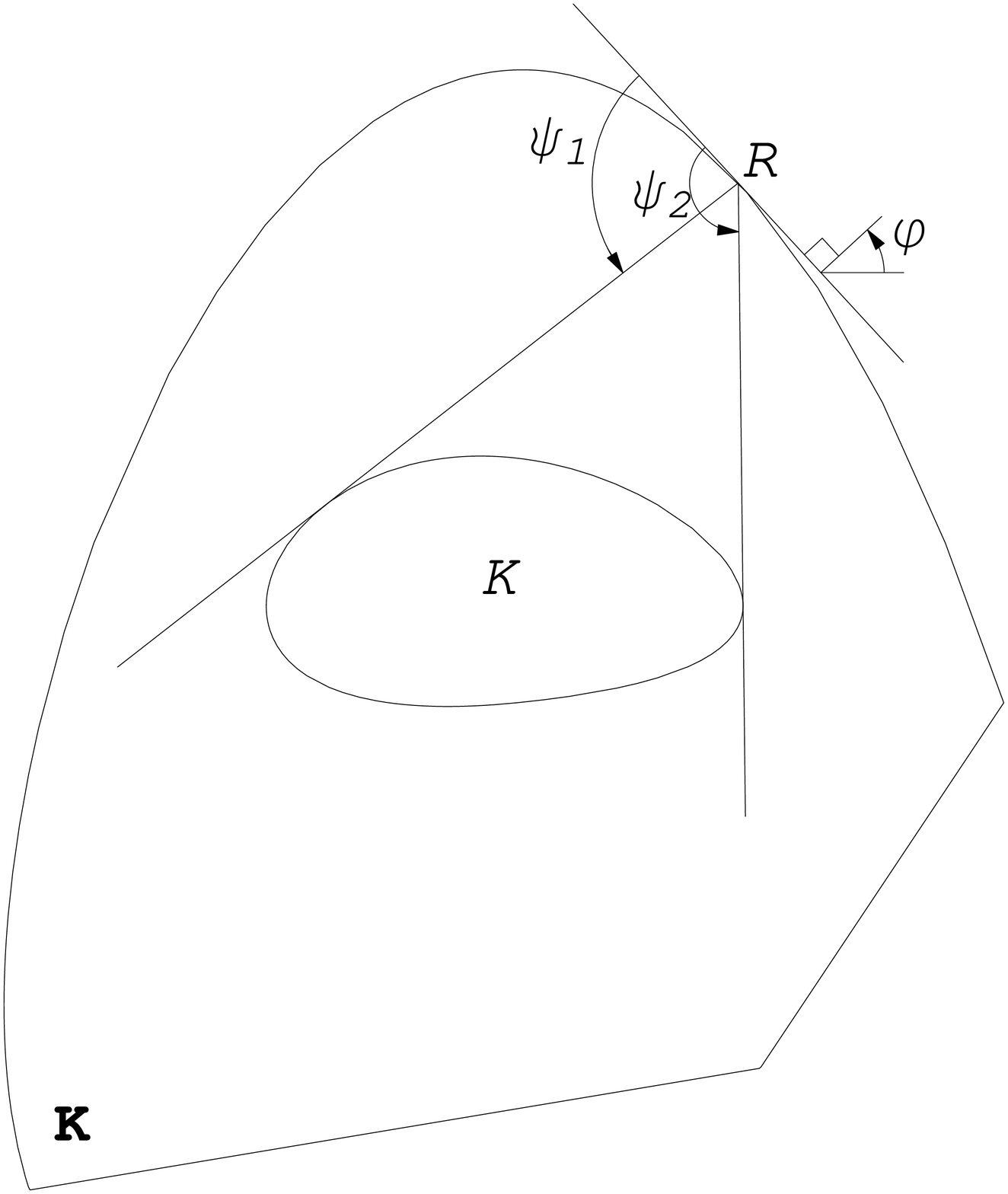}
\vskip-2.2in
\seteps{3.75in}{2in}{2.2in}{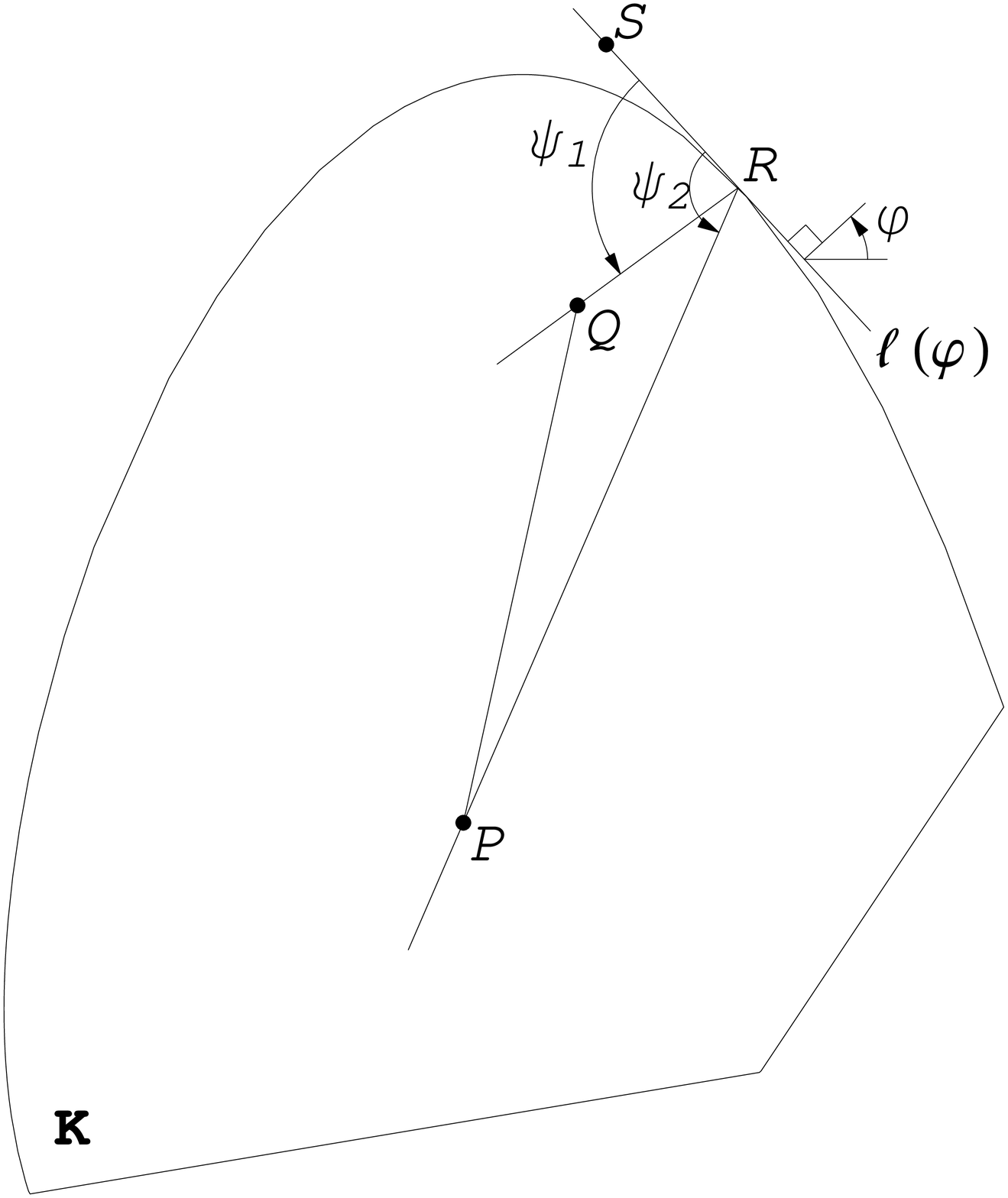}

\centerline{Figure~4. \hfil\hfil Figure~5. \hfil\hfil Figure~6.}

\smallskip

When $\bdy\ambientK$ is $C^2$, \thetag{1.10} can be reformulated as a Crofton-like formula. If $s$ is Euclidean arclength along $\bdy\ambientK$ and $\kappa = d\phi/ds$ is the Euclidean geodesic curvature of $\bdy\ambientK$, we have 
$$
\P = \frac12
\int_{\bdy\ambientK}\int_{\psi_1(s)}^{\psi_2(s)}\kappa\csc^2\psi\,d\psi\,ds
= \int_{\Cal L(K)}d\eta = \eta(\Cal L(K)),
$$
where $\Cal L(K)$ is the collection of lines passing through $K$ and $\eta$ is the measure on the oriented lines of the Hilbert plane induced by the density $\frac12\kappa\csc^2\psi\,d\psi\,ds$.

\subheading{Comparison of the Minkowski and Cauchy Formulas}
As mentioned earlier, in $\E2$ the two Minkowski forms \thetag{1.4} and the Cauchy form  and \thetag{1.3} are locally equivalent. Working informally we have
$$
\rho^2 \kg\,d\theta = r\kg\,ds = r\,d\omega,
$$
where the first equality follows from $r/\rho = \rho\,d\theta\!/\!ds$, a consequence of the similarity of the triangles in Figure~7.

\centereps{2.2in}{1.4in}{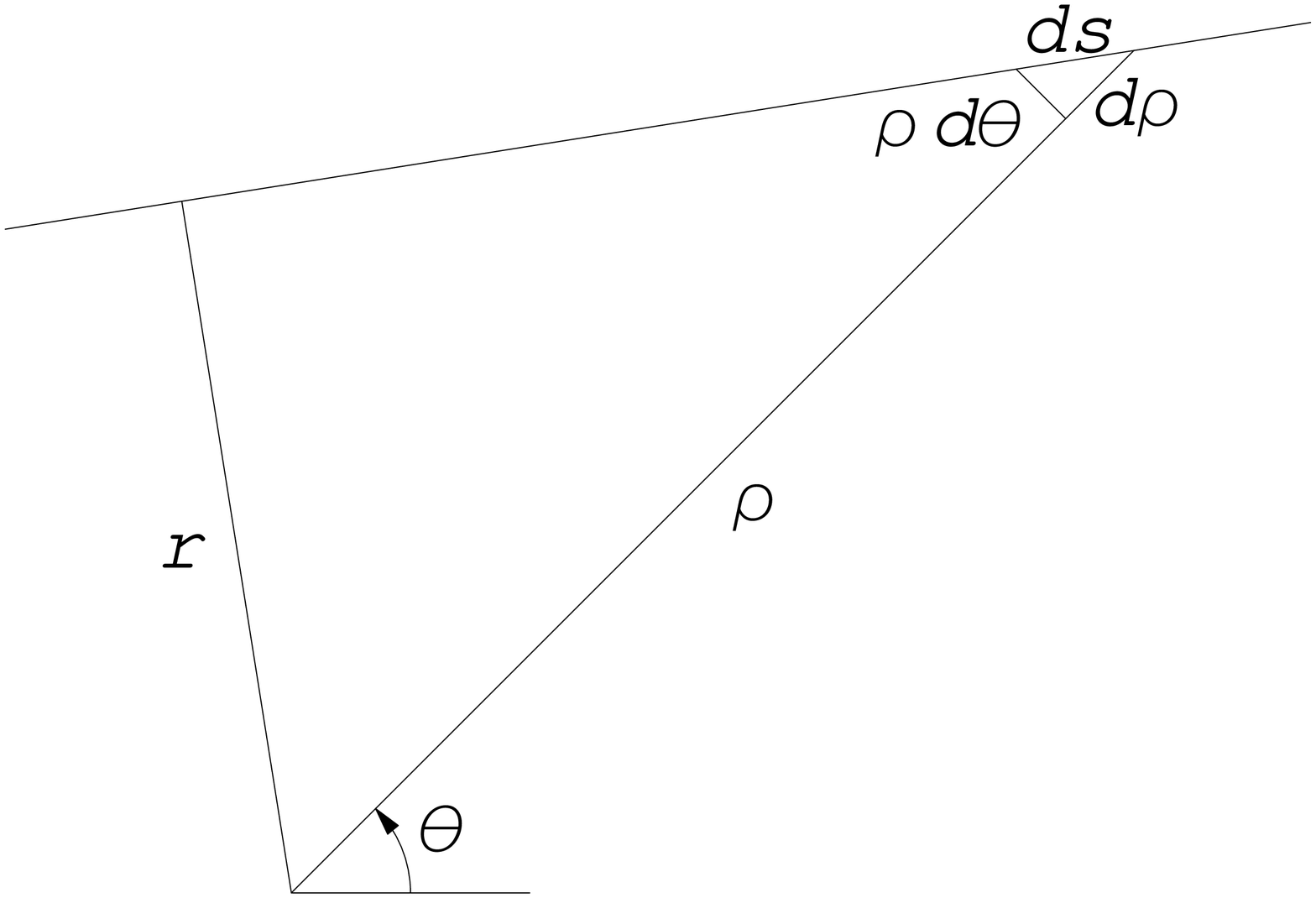}

\centerline{Figure~7.}

In contrast, \thetag{1.5} and \thetag{1.9}, are not locally equivalent in $\H2$ and $\S2$. For example, along a geodesic segment of $\bK$ the integrand of \thetag{1.9} vanishes (since $d\omega = 0$) as does the first part of the integrand of \thetag{1.5}, whereas the second part of the integrand of \thetag{1.5}, $k\frac{A(\rho)}{2\pi}\,ds$, does not. This comparison is most striking when the integrand of the Cauchy form is expressed in terms of $\theta$; it is very similar to the first part of the Minkowski integrand, and in particular is homogeneous in $\kg$. Indeed, the Cauchy formula \thetag{1.9} becomes (see \section5)
$$
\P = \int_\bK \(\frac{L(\rho)}{2\pi}\)^{\!2} \der\omega s\, d\theta 
  = \int_\bK \(\frac{L(\rho)}{2\pi}\)^{\!2} \kg\, \frac{2\pi - kA(x)}{2\pi - kA(r)}\, d\theta, \tag{1.12}
$$
where $x$ is the third side of the right triangle in Figure~3. While \thetag{1.5} and \thetag{1.12} are very similar and both reduce to the polar form of \thetag{1.4} when $k=0$, they have distinct integrands when $k\ne0$.

One step of the proof of \thetag{1.12} is the nice formula for the geodesic curvature of $\bK$, $\kg = \frac{2\pi - kA(r)}{2\pi - kA(x)}\der\omega s$, which becomes $\frac{\cos r}{\cos x}\der\omega s$ when $k=1$, $\frac{\cosh r}{\cosh x}\der\omega s$ when $k=-1$, and simply $d\omega/ds$ when $k=0$. 
In \section6 we derive other formulas relating $ds$, $d\phi$, $d\theta$, and $d\omega$, which can be viewed as relations between measures on $\bK$, and involve the geometry of $\bK$ and its location relative to the origin. 

We now set to work.

\heading 2. Minkowski Formulas \endheading

In this section we generalize the polar Minkowski formula \thetag{1.4} to $\S2$ and $\H2$. Let $M$ be $\E2$, $\S2$ or $\H2$ with constant curvature $k$. For $\rho\ge0$, let $L(\rho) = 2\pi\,\l(\rho)$ and $A(\rho) = 2\pi\,a(\rho)$ be, respectively, the circumference and area of a disk of radius $\rho$ in $M$. Formulas for these are in the appendix.

\proclaim{Theorem 2.1} Let $K$ be a compact region in $M$ with $C^2$ boundary. Let $\kg$ be the geodesic curvature of $\bK$. Choose a distinguished point $O$ in $M$ as origin, and let $(\rho,\theta)$ be polar coordinates relative to $O$. Then the perimeter of $K$ is
$$
\P = \int_\bK \l(\rho)^2\,\kg\,d\theta + k \int_\bK a(\rho)\,ds.
$$
\endproclaim

\demo{\smc Proof} 
When $k =0$ this is the problem Montes posed. We give a simple proof of
this special case as an illustration, and then generalize it. Let $T$
be the unit tangent vector for $\bK$ in the positive direction, and let $N = -T\Perp$, where $T\Perp$ the result of rotating $T$ counterclockwise by 90$^\circ$ (so $N$ is the outward normal to the boundary). Let $R$ be the position vector from the origin, that is, $R = \rho u$ is the Euler field, where $u = \grad\rho$ is the unit vector field pointing away from the origin. Assuming that $\bK$ does not pass through the origin, we have
$$
\der{}s\<R,T> = \<T,T> - \<R,\kg N> = 1 - \rho \kg \<u,N> 
= 1 - \rho^2\,\kg\,\der\theta s,
$$
where the final equality follows from
$\<u,N> = \<u\Perp,T> = \rho\,d\theta(T) = \rho\der\theta s$. Thus
Montes' integrand becomes $\rho^2\,\kg\,d\theta = \(1 - \der{}s\<R,T>\)\,ds$, and integrating around $\bK$ gives the result. Finally, note that $\bK$ can pass through the origin, since the form $\rho^2\,d\theta = -y\,dx + x\,dy$ is smooth there. 

 For the general case, let $T$, $N$, and $u = \grad \rho$ be as above. In
place of the Euler field $\rho u$ we set $R = \l(\rho)u$. We need to
compute $\der{}s\<R,T>$, which is accomplished by computing the
covariant derivatives $\del_Tu$ and $\del_TR$.

 We have $\del_uu = 0$, since the flow of $u$ is the geodesic flow away
from the origin. We also have $\del_{u\Perp}u = c(\rho)u\Perp$, where
$c(\rho)$ is the geodesic curvature of the circle of radius $\rho$. Assume
for the moment that $\bK$ doesn't pass through the origin $O$. On the
sphere we also require that $\bK$ not pass through the point $O'$
antipodal to $O$. Then $\del_Tu = \<T,u\Perp>c(\rho)u\Perp$ and
$$
\del_TR = \l'(\rho)\<T,u>u + \l(\rho)c(\rho)\<T,u\Perp>u\Perp.
$$
Now we use the formula $\l'(\rho) = \l(\rho)c(\rho) = 1 - ka(\rho)$. This can be proved either from a direct geometric analysis of how $\l(\rho)$, $a(\rho)$, $c(\rho)$, and $k$ are related, or by using the explicit formulas for these quantities (see \thetag{A.1} in the appendix). We get $\del_TR = \(1 - ka(\rho)\)T$. We now have
$$
\der{}s\<R,T> = \(1 - k a(\rho)\)\<T,T> - \<R,\kg N> = 1 - k a(\rho) - \l(\rho)\kg\<u,N>.
$$
Note that the quantity $\l(\rho)$ is the rate at which arc-length
accumulates with respect to $\theta$ when moving around the circle of
radius $\rho$. Thus $\<u,N> = \<u\Perp,T> = \l(\rho)\,d\theta(T) = \l(\rho)\der\theta s$, and so we have
$$
\der{}s\<R,T> = 1 - ka(\rho)-\l(\rho)^2 \kg \der\theta s.
$$
The integrand in the statement of the theorem thus becomes 
$\l(\rho)^2\,\kg\,d\theta\, +\, ka(\rho)\,ds\allowmathbreak = 
\allowmathbreak(1\nomathbreak-\nomathbreak\der{}s\<R,T>)\,ds$, and the result follows by integrating around
$\bK$. As before, the restriction that $\bK$ not pass through $O$ can be
lifted by noting that the form $\l(\rho)^2\,d\theta$ is
smooth there. Similarly, on the sphere the curve may pass through the
point $O'$. This concludes the proof of Theorem 2.1.
\enddemo

\remark{Remarks} 
\roster
\item This is really a theorem about the length of a closed curve as opposed to the perimeter of a set, that is, the proof goes through if $\bK$ is replaced by a $C^2$ closed curve.

\item The integrand may be expressed as a multiple of $d\theta$ by introducing the angle $\a$ from the radial vector $u$ to the outward normal $N$ (Figure~2). We have 
$\cos\a = \<u,N> = \<T,u\Perp> = \l(\rho)\der\theta s$, and so $ds = \l(\rho)\sec\a\,d\theta$. The integral then becomes
$$
\P = \int_\bK \(\l(\rho)^2\,\kg + k a(\rho)\l(\rho)\sec\a\)d\theta.
$$
This form will be particularly useful when we compare this Minkowski formula for perimeter with the Cauchy formulas. It is also useful when $K$ is star-like with respect to $O \in \text{int}(K)$, in which case $\theta$ can be used as a parameter for the curve and the integral.

\item The proof shows that the given integrand can be rewritten as 
$\big(1 + \der fs\big)ds$, where $f$ is some function on $\bK$, in fact $f = -\<R,T> = \l(\rho)\sin\a$. Thus the integrand can be written variously as
$$
\(1 - \der{}s\<R,T>\)\!ds = \(1 + \der{}s \(\l\(\rho\)\sin\a\)\)\!ds
 = \(\der s\theta + \der{}\theta \(\l\(\rho\)\sin\a\)\)\!d\theta.
$$
Similar computations yield formulas in higher dimensions and co-dimensions.

\item If $\bK$ is only piecewise $C^2$, then $0 = \int_\bK d\<R,T>$ must be interpreted in the Riemann-Stieltjes sense. On each $C^2$ piece $\bK_i$ of $\bK$ we get 
$
\int_{\bK_i}\! d\<R,T> = \int_{\bK_i}\! \(1 -ka\(\rho\)\)ds 
   - \int_{\bK_i}\!\l(\rho)^2\,\kg\,d\theta
$
as above. At vertex $v_i$ we get
$
\int_{v_i} d\<R,T> = \<R_i,\Delta T_i> = -\l(\rho_i)\Delta\sin\a_i
$.
After summing, the resulting formula for perimeter is
$$
\P = \int_\bK \l(\rho)^2\,\kg\,d\theta 
     + \sum_i \l(\rho_i)\Delta\sin\a_i + k \int_\bK a(\rho)\,ds.
$$

It is not our intent to explore the minimal regularity necessary for some version of \thetag{1.5} to hold, but there are two cases worth mentioning, when $\bK$ is piecewise $C^1$ or when $K$ is convex (in the latter case one-sided tangents to $\bK$ exist everywhere and are continuous except on a countable set). In either case one can make sense of the formula
$\P = -\int_\bK \<R,dT> + k \int_\bK a(\rho)\,ds$, which is in the spirit of the first part of the problem posed by Montes, namely that the perimeter of a convex body in $\E2$ is given by 
$-\int_\bK \<R,dT>$, appropriately interpreted \cite{M}. Here $-\<R,dT> = -\l(\rho)\<u,dT>$ is a measure on $\bK$ (signed if $\bK$ is not convex or if the origin is not contained in $K$), which can be defined using parallel translation or by approximating $\bK$ uniformly by polygons or $C^2$ curves, and which equals $\l(\rho)^2\,\kg\,d\theta$ at points where $\bK$ is $C^2$.

For Cauchy formulas, the vertices of a convex body, which are necessarily countable, do not cause such a technical problem. 
\endroster
\endremark

\heading 3. Cauchy Formulas for Hilbert Planes \endheading

In this section we derive the Cauchy formula \thetag{1.10} and the corresponding Crofton formula in two-dimensional Hilbert geometries. We assume the reader is familiar with the basics of projective geometry. See, for example, the book by Levy \cite{L}.

If $\ambientK$ is a bounded convex region in $\E{n}$ having nonempty interior, Hilbert \cite{H} pointed out that one may define a metric $h$ in the interior of $\ambientK$ by imitating the hyperbolic metric on the projective model of $\H{n}$. Thus, if $P$ and $Q$ are interior points and the line $\Line{PQ}$ strikes the boundary of $\ambientK$ at points $A$ and $B$ in the order $APQB$, define
$$
h(P,Q) = \frac12\log\(\frac{AQ}{AP}\,\frac{BP}{BQ}\) = -\frac12\log(P,Q;A,B),
$$
where $(P,Q;A,B)$ denotes the usual projective cross ratio (see Figure~4).
Here we only consider the case $n=2$ because the integral formulas developed in this section do not generalize to higher dimensional Hilbert geometries without very special assumptions about the body $\ambientK$. 

Next, define a projectively invariant pseudometric $d$ in the interior of a
given angle as follows. If $C$ is the vertex, choose points $A$ and $B$ on the
respective sides ${V}$ and ${W}$. Interior points $P^{\ast}$ and $Q^{\ast}$ determine rays ${T}$ and ${U}$ from $C$ which will cut segment $\seg{AB}$ at points $P$ and $Q$. We then define $d(P^{\ast},Q^{\ast}) = \frac12\log\(\frac{AQ}{AP}\frac{BP}{BQ}\)$. It
is clear that $d(P^{\ast},Q^{\ast})=0$ if and only if $C$, $P^{\ast}$, and
$Q^{\ast}$ are collinear.

As the pencil of lines through $C$ forms a one-dimensional projective space, the cross ratio $(T,U;V,W)$ of concurrent lines is canonically defined. A fundamental duality result says that the point cross ratio $(P,Q;A,B)$ and the line cross ratio $(T,U;V,W)$ are equal. There is the useful self dual notation $(P,Q;V,W)$ for this number, especially in light of the pretty formula
$$
(P,Q;V,W)=\frac{VP}{VQ}\,\frac{WQ}{WP}.
$$
Here the products $VP$, etc., are matrix products where $V=[v_1\;v_2\;v_{3}]$, and $P=[p_1\;p_2\;p_{3}]^{T}$ are homogeneous coordinates in $\Bbb P^2$ as discussed in \cite{L}, Chapter 3, for example. The embedding $(x,y) \mapsto [x\;y\;1]^{T}$ is the most natural way of embedding $\E2$ into $\Bbb P^2$. The associated dual embedding is given by $\{ax+by+c=0\} \mapsto[a\;b\;c]$.

The following discrete Cauchy formula for Hilbert metrics on planar polygonal
regions was stated without proof in \cite{A}. The proof offered here is the one in
mind for this earlier paper.

\proclaim{Lemma}
Let $\ambientK$ be bounded by a convex polygon in $\E2$ with
vertices $C_1,C_2,...,C_{n}$. If $h$ is the Hilbert metric on the interior
of $\ambientK$, then for any points $P$ and $Q$ in the interior,
$$
2h(P,Q)=\sum d_{i}(P,Q),
$$
where $d_{i}$ is the projective pseudometric, defined above, for the angle at
vertex $C_{i}$ .
\endproclaim

\demo{Proof}
The lemma is clearly valid if $P=Q$, so suppose that $P$ and $Q$ are distinct
points. Line $\Line{PQ}$ will strike the boundary of
$\ambientK$ at points $A$ and $B$ with order $APQB$. Next note that if $A$
or $B$ is a vertex $C_{i}$ , then $d_{i}(P,Q)=0$, and this vertex may be
ignored. The remaining vertices are separated into two nonempty sets by line
$\Line{PQ}$. Reindexing, suppose that the vertices on one side
of this line in order from point $A$ toward point $B$ are $C_1,...,C_{k}$.
We will show that $h(P,Q)=\sum_{i=1}^k d_{i}(P,Q)$. The
lemma then follows at once, as the same holds for the vertices on the opposite
side. 

We respectively denote the lines $\Line{AC_1}$ and $\Line{C_{k}B}$ by $V_{0}$ and $V_{k}$; otherwise denote $\Line{C_{i}C_{i+1}}$ by $V_{i}$ for $1 \le i \le k-1$. At vertex $C_i$ we have
$d_i(P,Q) = -\frac12\log\big((V_{i-1}P/V_{i-1}Q)(V_iQ/V_iP)\big)$. If
these equations are added, $1 \le i \le k$, one sees that
$$
\sum_{i=1}^{k}d_{i}(P,Q)
= \frac12\log\(\frac{V_{0}Q}{V_{0}P}\,\frac{V_{k}P}{V_{k}Q}\) = h(P,Q),
$$
as required. 
\enddemo

\remark{Remark}
As the proof of the lemma shows, the distance $h(P,Q)$ may be calculated by summing on one side only of the line $\Line{PQ}$. This property will be inherited by any integral formula based on the lemma.
\endremark

If $\ambientK$ is an arbitrary, planar, compact, convex body, for each angle $\phi$ ($\operatorname{mod}2\pi)$ there is a unique support line $\l(\phi)$ with outer normal $N(\phi)=(\cos\phi,\sin\phi)$. If $\l$ contains exactly one point of $\bdy\ambientK$, denote it by $R(\phi)$; such a point is termed
$\phi$-{\sl exposed.} If $\l(\phi)$ intersects $\bdy\ambientK$ in a line segment, an arbitrary point in this segment is chosen as $R(\phi)$; there are at most countably many such $\phi$ (see Figure~6).

If $P$ and $Q$ are points in the interior of $\ambientK$, the line $\l(\phi)$ together with the rays $\ray{R(\phi)P}$ and $\ray{R(\phi)Q}$
determine two counterclockwise angles. If $S$ is a point on $\l(\phi)$, lying on the positive side of $R$, these angles are $\angle SRP$ and $\angle SRQ$. The angles $\psi_1(\phi)$ and
$\psi_2(\phi)$ are defined as $\psi_1=\min(\angle SRP,\angle SRQ)$ and $\psi_2=\max(\angle SRP,\angle SRQ)$. It is clear that the functions $\psi_1$ and $\psi_2$ are continuous where $R(\phi)$ is well-defined, and that they have one-sided limits otherwise. It is also clear that 
$0<\psi_1\leq\psi_2<\pi$, and that $\psi_1=\psi_2$ if and only if $R$
is one of the two points $A,B$ where line $\Line{PQ}$ strikes
$\bdy\ambientK$. 

\proclaim{Theorem 3.1}
Let $\ambientK$ be a compact convex body in $\E2$, and let
$P$ and $Q$ be points in the interior of $\ambientK$. Then there is the
following integral formula for the Hilbert distance $h(P,Q)$%
$$
2h(P,Q) = \frac12\int_0^{2\pi}\(\cot\psi_1(\phi) -\cot\psi_2\(\phi\)\)\,d\phi.
$$
\endproclaim



\seteps{.5in}{2in}{2in}{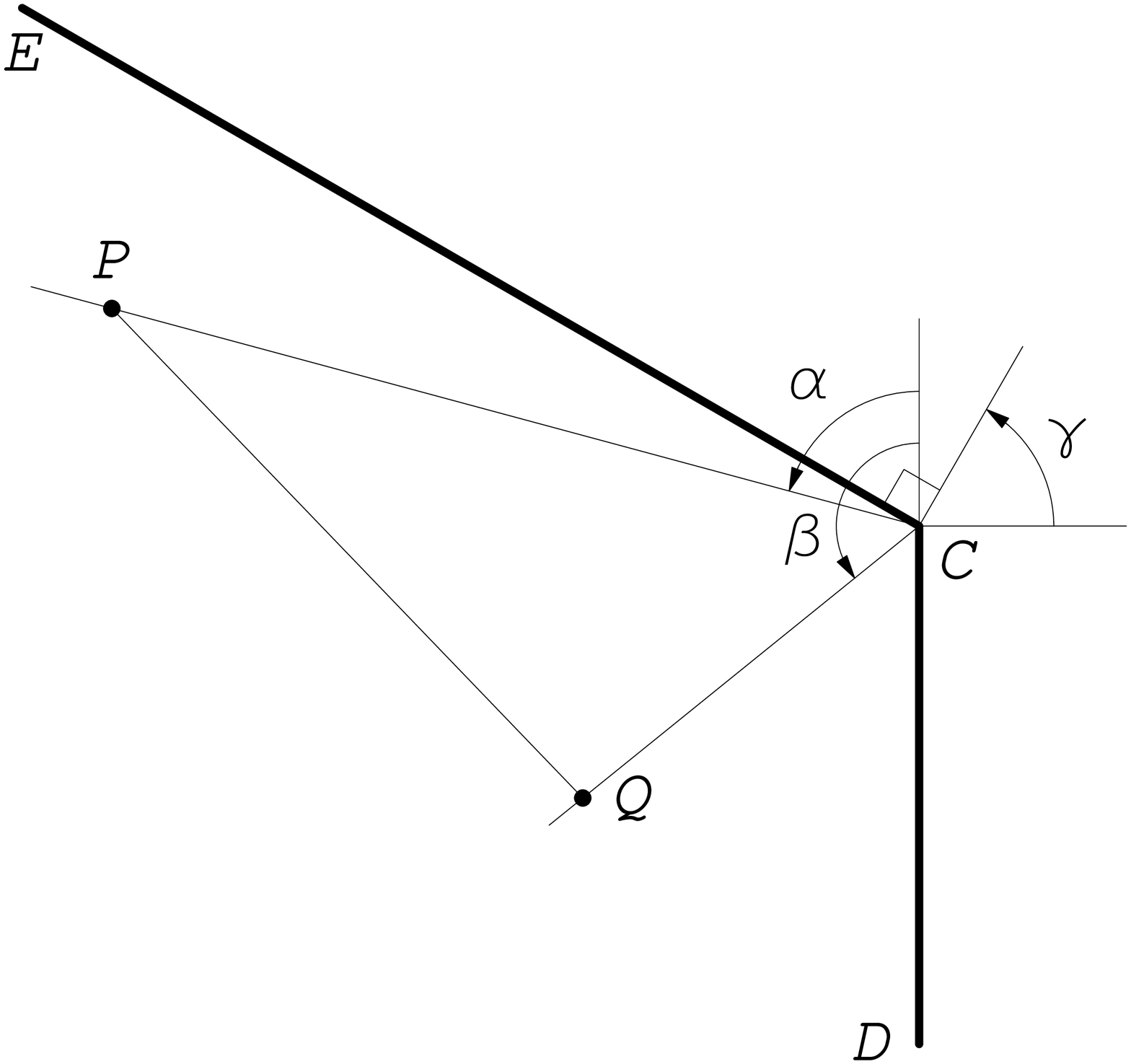}
\vskip-2.1in
\seteps{3in}{2in}{2in}{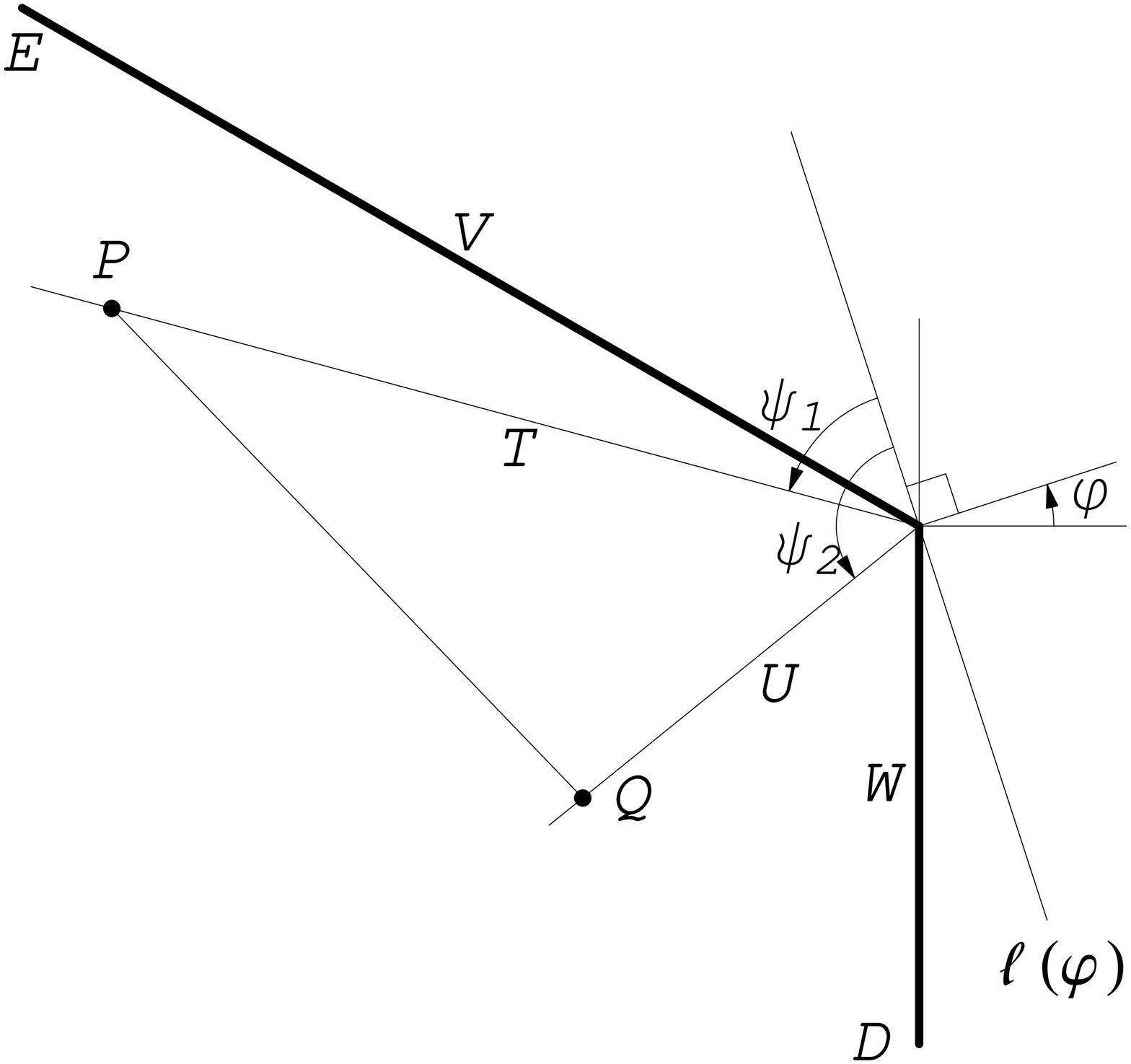}

\centerline{Figure~8. \hfil Figure~9.}
\bigskip

\demo{Proof}
The strategy is to first prove this formula when $\ambientK$ is a polygonal region. To achieve this we prove a version for the projective pseudometric $d$, defined above, on the interior of an angle $\angle DCE$ 
(see Figures~8 and~9). 
Suppose that this angle has measure $\pi-\gamma$ so that $\phi$ varies over an interval of length $\gamma$ at the vertex $C=R(\phi)$. Without loss of generality assume that $\phi$ varies from $0$ to $\gamma$. Furthermore, suppose that $\psi_1(0)=\alpha$ and $\psi_2(0)=\beta$, and so $\psi_1(\phi)=\alpha-\phi$ and $\psi_2(\phi)=\beta-\phi$. It is asserted that
$$
\gather
\int_0^\gamma\(\cot(\alpha-\phi)-\cot\(\beta-\phi\)\)\,d\phi
=\log\(\frac{\sin(\pi-\alpha)}{\sin(\pi-\beta)}
       \frac{\sin(\beta-\gamma)}{\sin(\alpha-\gamma)}\) \\
= -\log(T,U;V,W) = 2d(P,Q).
\endgather
$$
Here $W=\Line{CD}$ and $V=\Line{CE}$ are the sides of the angle, while $T=\Line{CP}$ and $U=\Line{CQ}$. The first equality above is a simple computation. The second equality follows from the
classical representation of the cross-ratio $(T,U;V,W)$ of four concurrent lines in the plane in terms of the sines of the various angles between the lines. \cite{K, pg.~196}

The result for an angle may now be combined with the lemma to assert that
$$
\frac12\int_0^{2\pi}\(\cot\psi_1(\phi) - \cot\psi_2\(\phi\)\)\,d\phi
= \sum d_i(P,Q) =2h(P,Q).
$$
Thus the theorem is valid when $\ambientK$ is a polygonal region.

For an arbitrary planar convex body $\ambientK$ it is possible to inscribe $n$-gons $\bdy\ambientK_{n}$ such that $\Delta(\ambientK,\ambientK_{n})\rightarrow0$; here $\Delta$ denotes Hausdorff distance.
For any $\phi$-exposed point $R(\phi)$ on $\bdy\ambientK$ it
will be true that $R_{n}(\phi)\rightarrow R(\phi)$. Also, $P$ and $Q$
are in the interior of $\ambientK_{n}$ for $n$ large, and it is clear that
$\psi_{ni}(\phi) \rightarrow \psi_{i}(\phi)$, $i=1,2$, at a $\phi$-exposed point. The exceptional values of $\phi$ where $R(\phi)$ is not
a $\phi$-exposed point form a countable set. Since $h_{n}(P,Q)\rightarrow
h(P,Q)$ and the angles $\psi_{ni}$ are bounded away from $0$ and $\pi$, the
theorem follows upon application of the bounded convergence theorem.
\enddemo

We now let $K$ be an arbitrary convex body in the interior of $\ambientK$. Given $\phi$, the angles $\psi_1$ and $\psi_2$ generalize in a natural manner (see Figure~5). The cone of rays from $R(\phi)$ through $K$ contains two extreme rays, which, with the support line $\l(\phi)$, determine two extreme
counterclockwise angles. Let $\psi_1(\phi)$ and $\psi_2(\phi)$ be the smaller and larger of these angles, respectively. A standard additivity argument in integral geometry leads at once to the following result.

\proclaim{Theorem 3.2}
If $K$ is a convex body in the interior of $\ambientK$, then the Hilbert perimeter of $K$ is 
$$
\P = \frac12\int_0^{2\pi}\(\cot\psi_1(\phi) - \cot\psi_2\(\phi\)\)\,d\phi,
$$
where $\psi_1(\phi)$ and $\psi_2(\phi)$ are the extreme angles subtended by $K$ with respect to $\l(\phi)$. 
\endproclaim

\demo{Proof} For a finite collection of segments inside $\ambientK$, say $\overline{P_{i}Q_{i}}$, it follows at once from Theorem 3.1 that $2\sum_{i}h(P_{i},Q_{i}) = \frac12\int_{0}^{2\pi}
\sum_{i}\(\cot\psi_{1}^{i}(\phi) - \cot\psi_{2}^{i}(\phi)\)\,d\phi$. If the  segments form the boundary of a convex polygon $K$ inside $\ambientK$, note that the sum inside the integral can be broken into two telescoping sums, each of which adds to $\cot\psi_{1}(\phi) - \cot\psi_{2}(\phi)$, where $\psi_{1}$, $\psi_{2}$ are the extreme angles for $K$. Thus the theorem is true for convex polygons, and the result for general convex bodies $K$ follows by a straightforward approximation argument.
\enddemo

With slightly more regularity of $\bdy\ambientK$, this can be expressed as a Crofton-type result, that is, as an integral with respect to a measure on the lines in the Hilbert geometry. As $d\phi$ defines a measure on $\bdy\ambientK$, the integral in Theorem~3.2 makes sense as an integral on $\bdy\ambientK$ except at the vertices, where the integrand is not well-defined. (However at each vertex, where $d\phi$ is $\Delta\phi$, the integrand can be replaced by its average value over the corresponding $\phi$-interval.) For simplicity, we assume that $\bdy\ambientK$ is $C^2$ and has positive Euclidean geodesic curvature $\kappa = d\phi/ds$, where $s$ is Euclidean arclength measured counterclockwise from some fixed point. Modifying our notation slightly, let $R(s)$ denote the point on $\bdy\ambientK$ with arclength parameter value $s$ and let $\l(s)$ be the unique support line passing through this point. The lines passing through $R(s)$ are parameterized by $0 < \psi < \pi$, where $\psi$ measures from $\l(s)$ as before. The {\sl oriented\/} lines of the Hilbert plane are then parameterized by $0 \le s < L$, $0 < \psi < \pi$, where $L$ is the Euclidean length $L$ of $\bdy\ambientK$. The set $\Cal L(K)$ of oriented lines passing through the convex body $K$ is parameterized by $0 \le s < L$, $\psi_1(s) \le \psi \le \psi_2(s)$, where $\psi_1(s)$ and $\psi_2(s)$ are the extreme angles as above. The formula in Theorem~3.2 becomes
$$
\P = \frac12\int_0^L\int_{\psi_1(s)}^{\psi_2(s)}\kappa\csc^2\psi\,d\psi\,ds.
\tag3.1
$$
We see that $\frac12\csc^2\kappa\psi\,d\psi\,ds$ is a density on the oriented lines, and so defines a measure $\eta$. Thus we have that $\P = \int_{\Cal L(K)} d\eta = \eta(\Cal L(K))$. More generally, for any $h$-rectifiable curve $\Cal{C}$, standard arguments \cite{S5, \section1.3.2} lead to the Crofton formula for the Hilbert length of $\Cal C$,

$$
L_h(\Cal C)=\frac12\int n(\l)\,d\eta(\l),
$$
where $n(\l)$ is the cardinality of $\l\cap\Cal{C}$.

Work related to Hilbert's fourth problem shows that within an open convex region $\ambientK$ of practically any 2-dimensional geodesic metric space there will be a unique measure $\eta$ on sets of geodesics such that length is given by a Crofton formula 
$d(P,Q)=\frac12 \eta\{\l\mid\nomathbreak\l\text{ cuts segment }\seg{PQ}\}$. As noted in \cite{A}, if $\ambientK$ is a polygonal region, the lemma gives a clear description of
$\eta$ as a sum of 1-dimensional projective measures on
projective angles. Such a projective angle measure is the
straightforward dual to the invariant measure on a hyperbolic
line. At first glance it might appear that for a smooth boundary $\bdy\ambientK$,
\thetag{3.1} gives a density function $d\eta$
for the lines of the geometry which can be integrated to obtain
$\eta(\Cal L)$ for a line set $\Cal L$. But the density
$\frac12\kappa\csc^2\psi\,d\psi\,ds$ is defined on
the {\sl oriented\/} lines so that each (unoriented) line is assigned
two such density values. However, the average of these densities will
be the desired density $d\eta$.

 In general a measure on the oriented lines naturally defines an
oriented distance function on the plane via an oriented Crofton
formula. The oriented line $(\l,\ray{u})$ is said to cut the oriented segment
$\ray{PQ}$ if (i) $\l\cap\seg{PQ}\neq 0$ and (ii) the counterclockwise
angle from $\ray{PQ}$ to $\ray{u}$ is less than $\pi$ radians. Define the
oriented distance $d(P,Q)$ to be half the measure of those oriented
lines which cut $\ray{PQ}$. For Hilbert planes, it is interesting to note that this oriented distance function defined by using the oriented line density 
$\frac12\kappa\csc^2\psi\,d\psi\,ds$ is symmetric, \ie, $d(P,Q)=d(Q,P)$, even though this function generally assigns two different values to a given
unoriented line. The proof of the lemma sheds light as to why this is true.

 Nonetheless, when $\ambientK$ is a circle, these oriented densities are
obviously equal on any given line. The circle is the most
important example because this construction leads to the projective
model for $\H2$, and here the measure $\eta$ will be invariant
under hyperbolic motions. This follows from the invariance of $h$
since $\eta$ is uniquely determined by $h$.

Hilbert geometries are among the earliest examples of what are now known as
Finsler spaces. Here the tangent spaces are Minkowski geometries (finite
dimensional real Banach spaces). In his 1903 treatment of Hilbert's fourth
problem Hamel \cite{Ha} gave an analysis of Hilbert planes from the viewpoint of the variational calculus. However, integral geometric methods generally have been more effective in dealing with the fourth and related problems. The paper of Chakerian \cite{Ch1} is the standard reference for integral geometry in Minkowski planes, and \cite{Ch2} describes more recent work in the same area.

\heading 4. Cauchy Formulas in $\H2$ \endheading

In this section we specialize the results of the previous section to $\H2$, and show that they imply the known Santal\'o-type Cauchy formula. In particular, we derive \thetag{1.6, 1.7, 1.8, 1.9}.




Let $\ambientK$ be the unit disk. The resulting Hilbert geometry, as developed in the previous section, is then the Beltrami-Klein model of $\H2$ with curvature $k=-1$. Let $K$ be a convex subset of $\H2$. For central angle $\phi$ the related quantities $R$, $\psi_1$, $\psi_2$, $w$, and $h$, defined in the introduction and the previous section, are shown in Figure~10 
(\cf Figures~2 and~5). 
Note that the $\phi$-normal is the line through the origin parallel to $\l(\phi)$.

One sees that $w(\phi) = \cot\psi_1(\phi) -\cot\psi_2(\phi)$ is the Euclidean length of the projection of $K$ from $R(\phi)$ onto the diametric line that is parallel to $\l(\phi)$. Theorem~3.2 then becomes the following.

\proclaim{Theorem 4.1 (Projective Cauchy formula)}
Let $\H2$ with curvature $k=-1$ be represented by the Beltrami-Klein model. Let $K$ be a convex body in $\H2$. Then
$$
\P = \frac12 \int_0^{2\pi} w \,d\phi,
$$
where $w$ is the Euclidean length of the projection of $K$ from $R$ onto the $\phi$-normal.
\endproclaim

\centereps{3in}{3in}{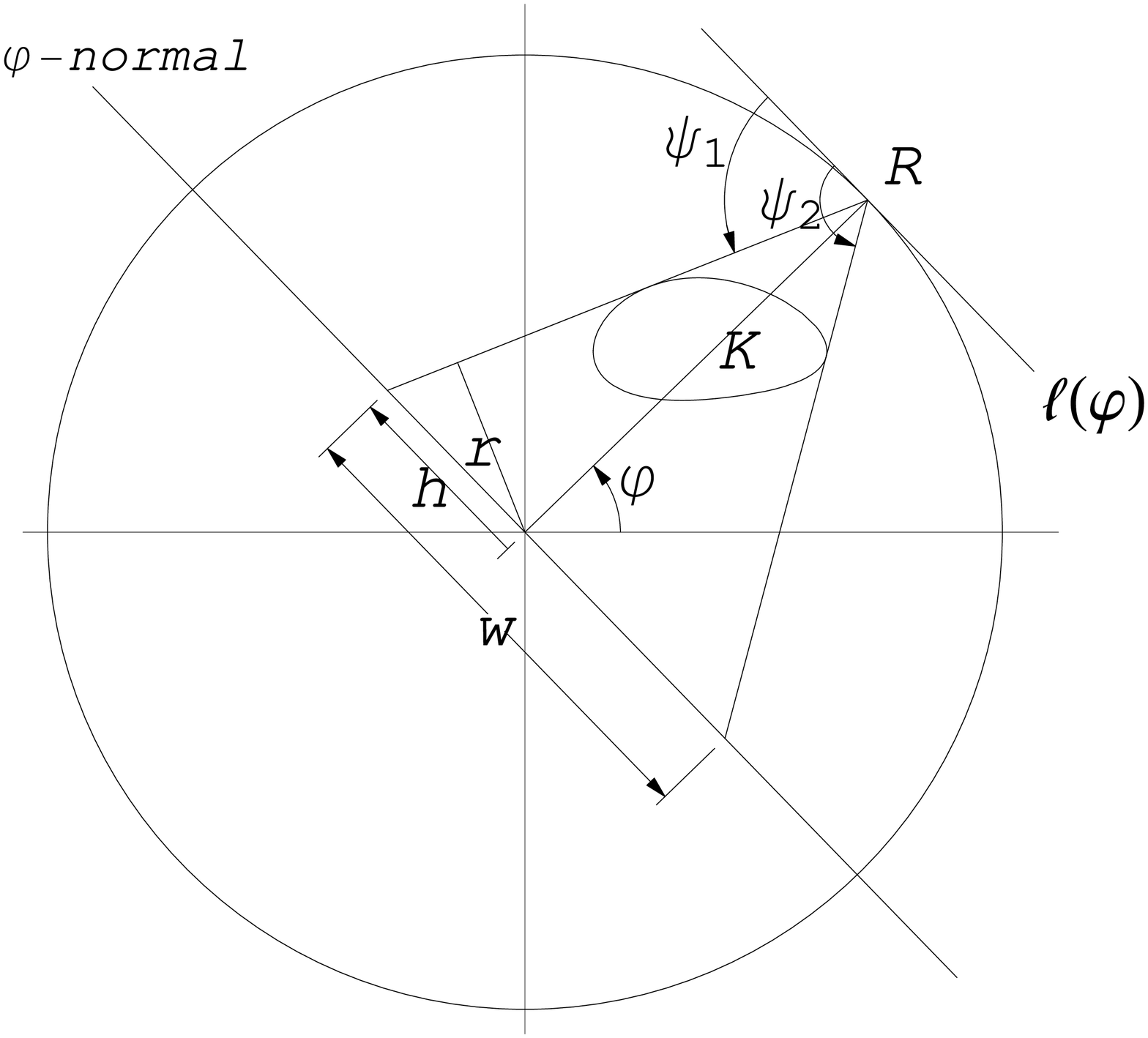}

\centerline{Figure~10.}


In addition, the symmetry of the circle allows a very useful one-sided version of the perimeter formula.

\proclaim{Theorem 4.2 (One-sided projective Cauchy formula)}
Let $\H2$ with curvature $k=-1$ be represented by the Beltrami-Klein model. Let $K$ be a convex body in $\H2$. Then
$$
\P = \int_0^{2\pi} h\,d\phi
$$
where $h$ is the signed Euclidean length of the projection of the right hand support line from $R$ onto the $\phi$-normal.
\endproclaim

\demo{Proof} Since $h(\phi) = \cot\psi_1(\phi)$, we must prove that $\P = \int_0^{2\pi} \cot\psi_1(\phi)\,d\phi$.

First, consider the case when $K$ is a segment $\seg{PQ}$, for which we need to prove that $2h(P,Q) = \int_0^{2\pi} \cot\psi_1(\phi)\,d\phi$, where $h$ denotes hyperbolic distance. By Theorem~3.1 we have 
$2h(P,Q) = \int_0^{2\pi} \(\cot\psi_1(\phi) - \cot\psi_2\(\phi\)\)\,d\phi$.
By rotational invariance, there is no loss of generality in assuming that $\seg{PQ}$ is perpendicular to the $x$-axis. The line $\Line{PQ}$ intersects $\bdy\ambientK$ in two points $A$ and $B$ determined by central angles $-\alpha$ and $\alpha$, respectively, for some $\alpha$ in $(0,\pi)$. By the remark following the lemma in the previous section, $h(P,Q)$ can be computed by integrating on just one side of $\seg{AB}$, \ie, 
$$
2h(P,Q) 
= \int_{-\alpha}^\alpha \(\cot\psi_1(\phi) - \cot\psi_2\(\phi\)\)\,d\phi
= \int_\alpha^{2\pi-\alpha} \(\cot\psi_1(\phi) - \cot\psi_2\(\phi\)\)\,d\phi.
$$
We further specialize by assuming $P=(x,0)$ and $Q=(x,y)$, $y>0$. Observe that, by symmetry, $\int_{-\alpha}^{\alpha}\cot\psi_2(\phi)\,d\phi
=\int_{\alpha}^{2\pi-\alpha}\cot\psi_1(\phi)\,d\phi=0$, as $\psi_2(-\phi) = -\psi_2(\phi)$ for $-\alpha\leq\phi\leq\alpha$ and $\psi_1(2\pi-\phi) = -\psi_1(\phi)$ for $\alpha<\phi<2\pi-\alpha$. The expression for $h(P,Q)$ reduces to 
$2h(P,Q) = \int_0^{2\pi} \cot\psi_1(\phi)\,d\phi 
= \int_0^{2\pi} \cot\psi_2(\phi)\,d\phi$, which gives the result for this special segment. If $P=(x,0)$ and $Q=(x,y)$, $y<0$, the same argument,
{\sl mutatis mutandis}, gives the result in this case. Finally, we note that any segment perpendicular to the $x$-axis is the set union or the set difference of two such special segments, and so the result holds for an arbitrary segment. The result extends to general perimeters of convex sets by additivity.
\enddemo

We now turn to the intrinsic realization of our Projective
Cauchy formula. Let $\H2$ be the hyperbolic plane with curvature $k=-1$, and let $K \subset \H2$ be a convex body. Choose some point $O \in \H2$ to be the origin. As above and in the introduction, given a central angle $\phi$ at $O$ measured from some fixed direction, let $R$ be the corresponding point at infinity, and consider the right-hand support line of $K$ passing through $R$ as viewed from $R$. Let $r$ be the signed hyperbolic distance from $O$ to the support line, let $\rho$ be the hyperbolic distance from $O$ to a point $P$ of contact of the support line with $K$, and let $\alpha$ be the hyperbolic angle from the radial direction to the outward normal at $P$ (see Figure~3). We take the sign of $r$ to be that of $\cos\alpha$; in fact, $\sinh r = \sinh\rho\,\cos\alpha$ by the law of sines. We get the following intrinsic version of Theorem~4.2.

\proclaim{Theorem 4.3} In $\H2$ with curvature $k=-1$, the perimeter of a convex body $K$ is given by
$$
\P = \int_0^{2\pi} \sinh r\,d\phi = \int_0^{2\pi} \sinh\rho\,\cos\a\,d\phi.
$$
\endproclaim

\demo{Proof} Representing $\H2$ by the Beltrami-Klein model with $O$ at the center of the model, from the previous theorem we have
$$
\P = \int_0^{2\pi} h\,d\phi = \int_0^{2\pi} \cot\psi_1(\phi)\,d\phi,
$$
where $h$ and $\psi_1$ are defined as above (see Figure~10). Letting $r_E$ be the Euclidean distance from $O$ to the support line, we have $r = \frac12\log\(\frac{1+r_E}{1-r_E}\) = \tanh^{-1}r_E$, and so $h = \cot\psi_1 = r_E/\sqrt{1 - r_E^2} = \sinh r$, which relates the extrinsic quantities $h$ and $\psi_1$ to the intrinsic quantity $r$, and completes the proof.
\enddemo

We note in passing a second intrinsic interpretation of Theorem~4.2. Let $c_1$ and $c_2$ be oriented circles in the extended complex plane (as usual, we include Euclidean lines in this family). Associated with these circles is the real number $(c_1,c_2)$, their signed inversive product, which is invariant under the group of M\"obius transformations. Considering either of the conformal models of $\H2$, the Poincar\'e disk or upper half-plane, it makes sense to talk about the signed inversive product of two oriented lines $\l_1$ and $\l_2$ in $\H2$. If $\l_1$ and $\l_2$ intersect, then $(\l_1,\l_2) = \pm\cos\theta$, where $\theta$ is the angle between them. If they don't intersect, then $(\l_1,\l_2) = \pm\cosh d$, where $d$ is the distance between them. In each case the sign is determined by the orientations of the lines. It can be shown that the signed Euclidean length $h$ in Theorem~4.2 equals $(\Line{OS},\Line{PR})$, where $\Line{OS}$ and $\Line{PR}$ are, respectively, the $\phi$-normal and the right-hand support line with appropriate orientations (see Figure~2). Many of the details concerning inversive products are contained in \cite{B}. Those remaining are left to the interested reader.

We now deduce the known Cauchy formula for $\H2$ \cite{S4}. Let $\omega$ be the central angle for the point on the support line closest to $O$, as in Figures~2 and~3. 

\proclaim{Corollary} In $\H2$ with curvature $k=-1$, the perimeter of a convex body $K$ is given by
$$
\P = \int_0^{2\pi} \sinh r\,d\omega.
$$ 
\endproclaim

\centereps{3in}{1.5in}{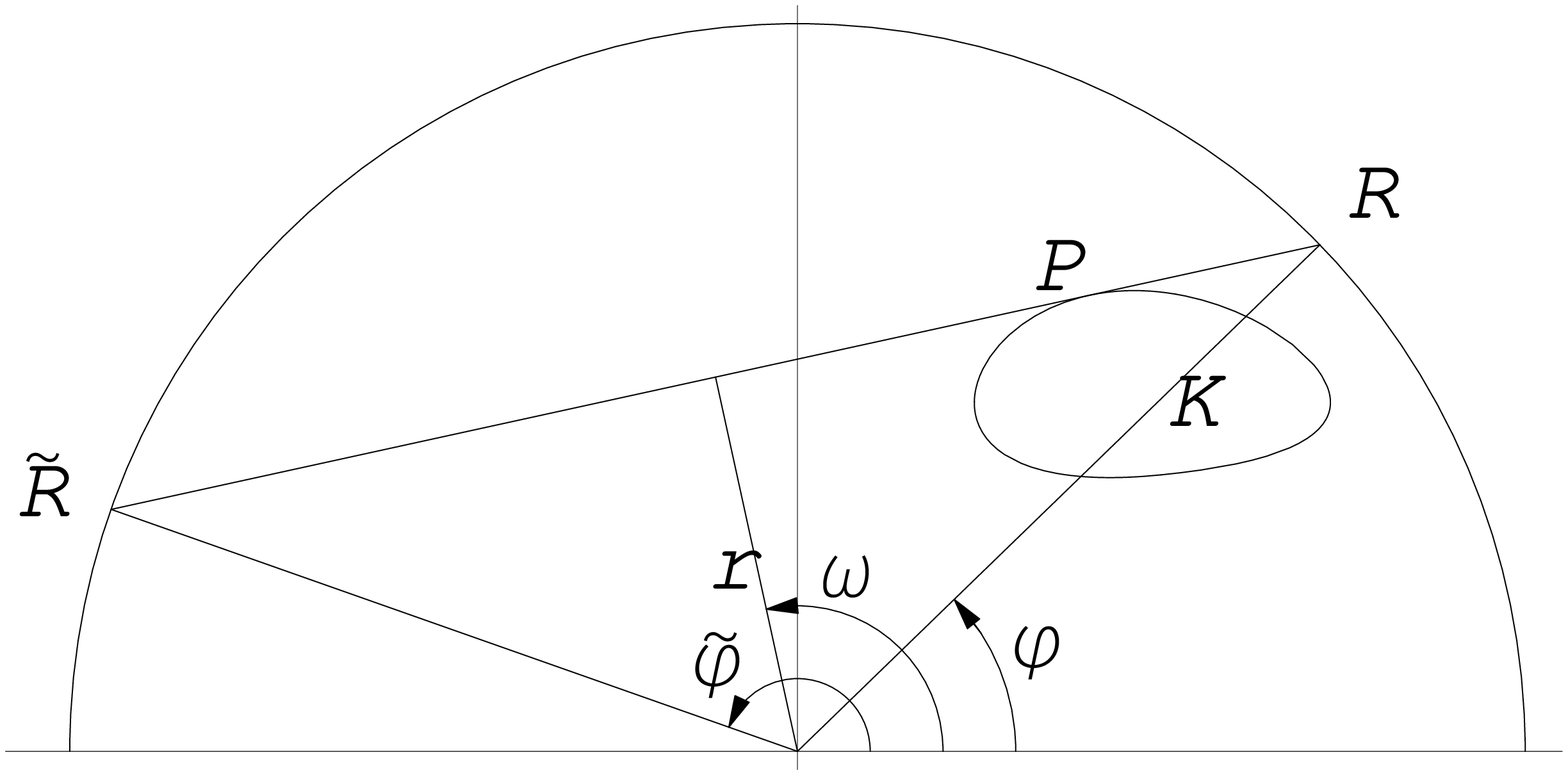}

\centerline{Figure~11.}

\demo{Proof}From Theorem~4.3 we have $\P = \int_0^{2\pi}\sinh r\,d\phi$. Given $\phi$, which determines $R$ and the right-hand support line $\Line{PR}$ as above, let $\tilde R$ be the opposite ideal endpoint of $\Line{PR}$, with central angle $\tilde\phi$ (see Figure~11). Then $\Line{PR}$ is the left-hand support line of $K$ as viewed from $\tilde R$, and the same argument leading up to Theorem~4.3 implies that
$\P = \int_0^{2\pi}\sinh r\,d\tilde\phi$. The result follows since $\omega = (\phi + \tilde\phi)/2$.
\enddemo



Note that the support line in the corollary can be parameterized by $\omega$ as opposed to $\phi$ as in Figure~3, in which case the corresponding construction works in $\E2$ and $\S2$ as well. (In $\S2$ one needs to assume that $K$ is contained in the open hemisphere centered at the origin $O$.) The corresponding formula in $\E2$ is \thetag{1.3}; in $\S2$ with curvature $k=1$ it is $\P = \int_0^{2\pi} \sin r\,d\omega$ \cite{S3}. In all three geometries with arbitrary constant curvature, the perimeter is given by the unified Cauchy formula
$$
\P = \int_0^{2\pi} \frac{L(r)}{2\pi}\,d\omega,
$$
where $L(r)$ is the circumference of the circle of radius $r$.

\bigskip




\heading 5. Comparison of the Cauchy and Minkowski Forms \endheading

In this section we compare the Cauchy and Minkowski forms
$$
\P = \int_0^{2\pi} \frac{L(r)}{2\pi}\,d\omega \qqand
\P = \int_\bK \left(\frac{L(\rho)}{2\pi}\right)^2\! \kg \,d\theta 
    + k\int_\bK \frac{A(\rho)}{2\pi}\,ds,
$$
where $L(r)$ and $A(r)$ are the circumference and area, respectively, of a circle of radius $r$. As noted in the introduction, the integrands of these are necessarily inequivalent in $\S2$ and $\H2$ (in contrast to $\E2$), since, in particular, along a geodesic segment of $\bK$ we have $\frac{L(r)}{2\pi}\,d\omega = 0$ and $\left(\frac{L(\rho)}{2\pi}\right)^2\! \kg = 0$, whereas $k\frac{A(\rho)}{2\pi} \ne 0$. Additional insight is obtained by expressing the Cauchy formula in terms of the polar angle $\theta$. 

The situation in all three geometries is given in Figure~3, in which we assume $\bK$ is $C^2$. In particular, $(\rho,\theta)$ are the polar coordinates of the point $P\in\bK$, $r$ is the distance from the origin to the unique support line of $K$ passing through $P$, $\omega$ is the central angle of the support line's closest point to the origin, and $\alpha$ is the angle from the radial vector $u$ to the outward normal $N$ to $\bK$ at $P$. Let $x$ be the length of the third side of the right triangle, opposite the angle $\omega - \theta$, and let $\beta = \pi/2 - \alpha$ be the angle opposite the side $r$. The proofs in this section make extensive use of identities \thetag{A.2} and \thetag{A.3} for this triangle (in the appendix). Note that $\omega$ can be taken to be continuous when the support line passes through the origin. Also note that $r$, $x$, and $\beta$ are signed quantities: the signs of $r$ and $x$ are taken to be the same as those of $\beta$ and $\omega-\theta$, respectively.

We begin with a pleasant formula for the curvature of $\bK$ in terms of the angle $\omega$. As in \section2 we simplify the notation by letting  $\l(r) = L(r)/2\pi$ and $a(r) = A(r)/2\pi$. We also let $c(r)$ be the geodesic curvature of the circle of radius $r$ (see the appendix). 

\proclaim{Theorem 5.1} The geodesic curvature of $\bK$ is given by 
$\dsize \kg = \frac{\l(r)c(r)}{\l(x)c(x)}\der\omega s 
  = \frac{1 - ka(r)}{1 - ka(x)}\der\omega s$.
\endproclaim

\demo{Proof} From \thetag{A.3} we have 
$\l(r) = \l(\rho)\cos\a = \l(\rho)\<u,N> = \<R,N>$, where $R = \l(\rho)u$. Recall from the proof of Theorem~2.1 that the covariant derivative $\frac{DR}{ds} = \del_T R$ is a multiple of $T$, where $T$ is the unit tangent to $\bK$ taken in the positive direction. From \thetag{A.1} we have $\l'(r) = \l(r)c(r)$. From the lemma below we have $dr/d\omega = -\l(r)c(r)/c(x)$. Using these observations, we differentiate both sides of $\l(r) = \<R,N>$ with respect to $s$, obtaining
$$
-\frac1{c(x)}\(\l(r)c\(r\)\)^2 \der\omega s = \<R,\kg T> = \kg\l(\rho)\<u,T>
 = -\kg\l(\rho)\cos\b = -\kg\l(\rho) \frac{c(\rho)}{c(x)},
$$
where the last equality uses \thetag{A.3}. Rearranging and applying \thetag{A.2} yields
$$
\l(r)c(r) \der\omega s = \kg \frac{\l(\rho)c(\rho)}{\l(r)c(r)} = \kg\l(x)c(x).
$$
Finally, $\l(r)c(r) = 1 - ka(r)$, as noted in the appendix.
\enddemo

\proclaim{Lemma} $\dsize \der r\omega = -\frac{\l(r)c(r)}{c(x)}$.
\endproclaim
\demo{Proof} First note that from \thetag{A.1} we have $c'(r) = -1/\l(r)^2$, which is a nice analogue of $\der{}r (1/r) = -1/r^2$. From \thetag{A.3} we have 
$\cos(\omega - \theta) = c(\rho)/c(r)$. Differentiating this with respect to $\omega$ and rearranging we get
$$
\frac1{\l(r)^2} \cos(\omega-\theta)\der r\omega + c(r)\sin(\omega-\theta)
=\frac1{\l(\rho)^2}\der\rho\omega + c(r)\sin(\omega-\theta)\der\theta\omega.
$$
Using $\sin\a = \cos\b = \l(r)c(r)\sin(\omega-\theta)$ from \thetag{A.3} the right hand side becomes
$$
\align
\frac1{\l(\rho)^2}\der\rho\omega & + \frac1{\l(r)}\der\theta\omega\sin\a
=\frac1{\l(\rho)\l(r)}
\(\frac{\l(r)}{\l(\rho)}\der\rho\omega + \,\l(\rho)\der\theta\omega\sin\a\)\\
&=\frac1{\l(\rho)\l(r)}
   \(\der\rho\omega\cos\a + \,\l(\rho)\der\theta\omega\sin\a\)
 =\frac1{\l(\rho)\l(r)} \<T,N>\der s\omega = 0.
\endalign
$$
Thus
$$
\cos(\omega-\theta)\der r\omega + c(r)\l(r)^2 \sin(\omega-\theta) = 0.
$$
Once again from \thetag{A.3} we have $\cos(\omega-\theta) = c(\rho)/c(r)$ and 
$\l(r)c(r)\sin(\omega-\theta) = \cos\b = c(\rho)/c(x)$, which give the result.
\enddemo

Finally we have the main result of this section.

\proclaim{Theorem 5.2} For a smooth convex body $K$ with $C^2$ boundary in $\E2$, $\S2$, or $\H2$ we have
$$
\dsize \P = \int_\bK \l(\rho)^2 \der\omega s \,d\theta 
   = \int_\bK \kg \l(\rho)^2 \frac{\l(x)c(x)}{\l(r)c(r)} \,d\theta
   = \int_\bK \kg \l(\rho)^2 \frac{1 - ka(x)}{1 - ka(r)} \,d\theta.
$$
\endproclaim
\demo{Proof} The second equality follows from Theorem~5.1. The first equality is a simple modification of the Cauchy formula \thetag{1.9}, $\P = \int_0^{2\pi} \l(r)\,d\omega$. We have $\l(r)/\l(\rho) = \sin\b =  \l(\rho)\,d\theta/ds$, where the first equality is from \thetag{A.3}. Thus
$$
\P = \int_0^{2\pi} \l(r)\,d\omega
  = \int_0^{2\pi} \l(\rho)^2 \der\theta s\,d\omega
  = \int_\bK \l(\rho)^2 \der\omega s\,d\theta.
$$
\enddemo
Note that when $k=0$ the proof of Theorem~5.2 reduces to the demonstration of the equivalence given in \thetag{1.11} of the Minkowski and Cauchy forms \thetag{1.3, 1.4} in the Euclidean case. The comparison of Theorems~2.1 and 5.2 shows that the Cauchy and Minkowski forms are genuinely inequivalent when $k\ne0$.

\heading 6. Comparison of $d\theta$, $d\omega$, $d\phi$, and $ds$ on $\bK$. \endheading

The angles $\theta$, $\omega$, and $\phi$ parameterize the circle of directions in $\E2$, $\H2$, and $\S2$ as viewed from the origin. In $\E2$ and $\H2$ the circle of directions is identified with the points at infinity, as we have been doing. In $\S2$ it is identified with the great circle centered at the origin in the obvious way. These parameters are related by the convex body $K$. (On $\S2$ we assume, as usual, that $K$ is contained in the open hemisphere centered at the origin. The construction relating $\phi$ to the other two angles is similar to that in $\E2$ and $\H2$. Some spherical trigonometry reveals that $\omega-\phi \equiv \pi/2$ on $\S2$, just as in $\E2$. Details are left to the reader.) The correspondence between $\omega$ and $\phi$ is always bijective. The relationship between $\theta$ and the other two angles is more complicated, but is bijective when the origin is in the interior of $K$ and $\bK$ is smooth and strictly convex.

The corresponding one-forms $d\theta$, $d\omega$, and $d\phi$ can be viewed as measures on $\bK$. (Actually, $d\theta$ is a measure when the origin is in the interior of $K$ and is a signed measure when the origin is in the exterior of $K$.) The ways these are related to each other and to the measure $ds$ involve the geometry of $\bK$ and its location relative to the origin. For simplicity, we assume that $\bK$ is $C^2$, but we do not assume that the origin is in the interior of $K$.

We have used $d\theta/ds = \l(r)/\l(\rho)^2$ several times already. We computed $\der\omega s = \kg \frac{\l(x)c(x)}{\l(r)c(r)}$ in Theorem~5.1. Combining these we have
$$
\der\omega\theta = \kg \frac{\l(x)c(x)}{\l(r)c(r)}\,\frac{\l(\rho)^2}{\l(r)},
$$
which was indirectly used in the proof of Theorem~5.2. We note that these hold in all three geometries.

The angle $\phi$ plays its most important role in hyperbolic geometry. Note that $\omega-\phi$ is the angle of parallelism for the length $r$. From \thetag{A.4} we have
$$
\cos(\omega-\phi) = \frac{\sqrt{-k}}{c(r)} \qqand 
\cot(\omega-\phi) = \sqrt{-k}\l(r).
$$
Differentiating $c(r)\cos(\omega-\phi) = \sqrt{-k}$ with respect to $\omega$ and using $c'(r) = -1/\l(r)^2$ and $dr/d\omega = -\l(r)c(r)/c(x)$ (from the lemma in \section5), we get
$$
\frac1{\l(r)^2}\frac{\l(r)c(r)}{c(x)}\cos(\omega-\phi) 
  - c(r)\sin(\omega-\phi)\(1 - \der\phi\omega\) = 0.
$$
After rearranging this becomes
$$
\der\phi\omega = 1 - \frac{\sqrt{-k}}{c(x)} = 1 - \tanh(\sqrt{-k}x).
$$
Combining this with $d\omega/d\theta$ we have
$$
\der\phi\theta 
= \(1 - \frac{\sqrt{-k}}{c(x)}\)\kg 
  \frac{\l(x)c(x)}{\l(r)c(r)}\,\frac{\l(\rho)^2}{\l(r)}. \tag6.1
$$

It is also instructive to compute $d\phi/d\theta$ directly, which we do in the Poincar\'e model of $\H2$ with curvature $k=-1$ (Figure~12).

\centereps{4in}{2in}{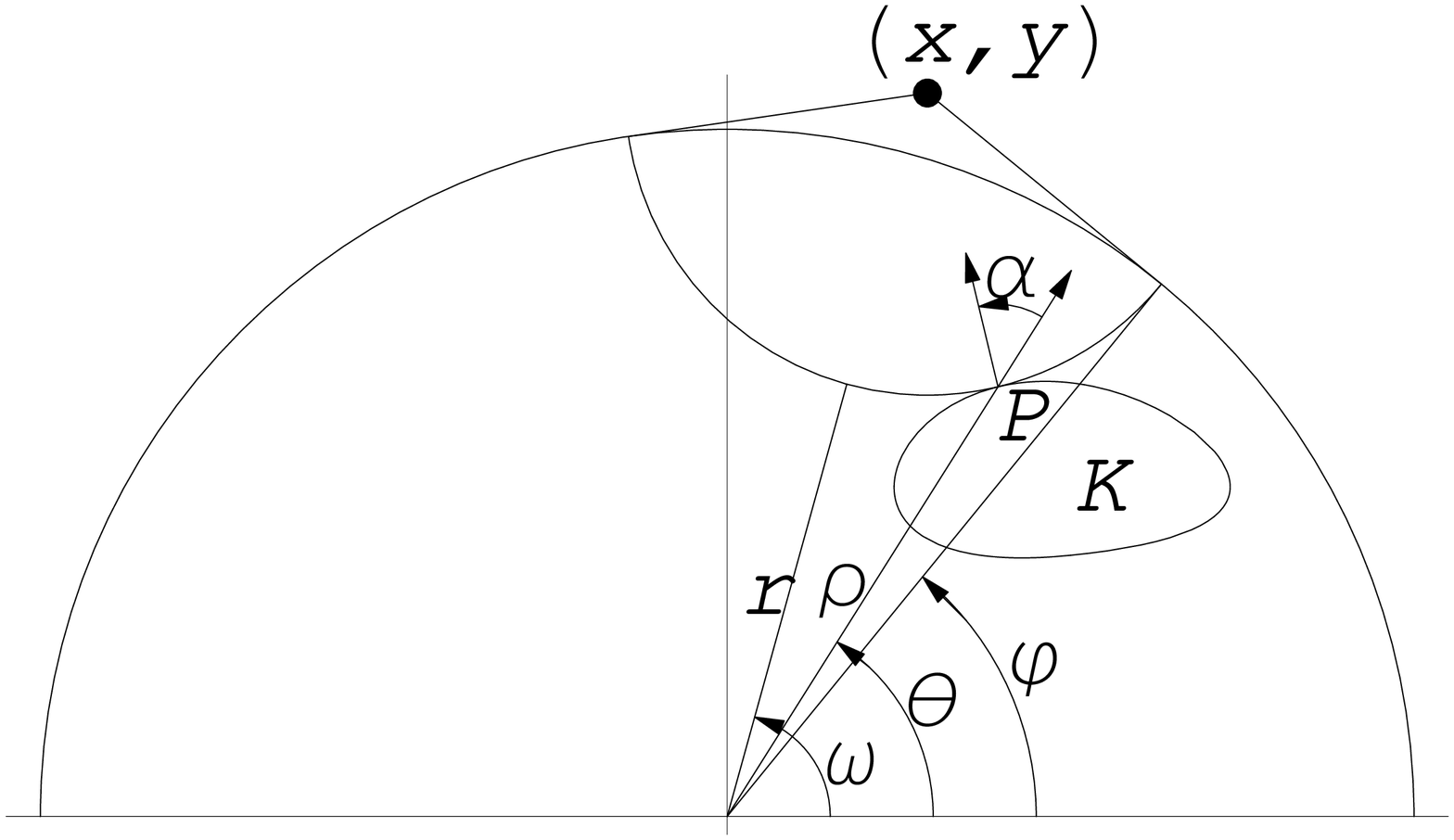}

\centerline{Figure~12.}
\medskip

The support geodesic is (assuming the generic case when it doesn't pass through the origin) the Euclidean circle of radius $\sqrt{x^2 + y^2 -1}$ centered at $(x,y)$, where
$$
x = \frac{\cos\theta\cosh\rho - \sin\theta\tan\a}{\sinh\rho} \qqand
y = \frac{\sin\theta\cosh\rho + \cos\theta\tan\a}{\sinh\rho}.
$$
Noting that $x\cos\phi + y\sin\phi = 1$, we differentiate with respect to arclength along $\bK$, obtaining
$$
\der\phi s (y\cos\phi - x\sin\phi) + \der xs\cos\phi + \der ys\sin\phi = 0. \tag6.2
$$
We seek to express everything in this equation in terms of $r$, $\rho$, and $\sin\a$. 

From the components of the unit tangent vector we have
$$
\der\rho s = -\sin\a \qqand 
\der\theta s = \frac{\cos\a}{\sinh\rho} = \frac{\sinh r}{\sinh^2 \rho},
$$
where $\cos\a = \sin\b = \sinh r/\sinh\rho$ by \thetag{A.3}. The geodesic curvature is given by
$$
\kg = \cos\a \frac{\cosh\rho}{\sinh\rho} + \der\a s.
$$

Next, we use $\sin(\phi-\theta) = 
\sin(\omega-\theta)\cos(\omega-\phi) - \cos(\omega-\theta)\sin(\omega-\phi)$ and the corresponding formula for $\cos(\phi-\theta)$. Since $\omega-\phi$ is the angle of parallelism for the distance $r$, we have (see the appendix or \cite{B})
$$
\sin(\omega-\phi) = \frac1{\cosh r} \qqand 
\cos(\omega-\phi) = \tanh r.
$$
From \thetag{A.3} we have
$$
\sin(\omega-\theta) = \frac{\cos\b}{\cosh r} = \frac{\sin\a}{\cosh r} \qqand
\cos(\omega-\theta) = \frac{\tanh r}{\tanh\rho}.
$$
Using these we get
$$
\sin(\phi-\theta) = \frac{\sin\a\sinh r - \cos\a\cosh\rho}{\cosh^2 r} \qand
\cos(\phi-\theta) = \frac{\cos\a\cosh\rho\sinh r + \sin\a}{\cosh^2 r}.
$$

Once the derivatives are taken, we lose no generality by assuming $\theta = 0$. We get
$$
\der xs = 0 \qqand \der ys = \frac{\kg}{\cos^2\a \sinh\rho}
= \kg \frac{\sinh\rho}{\sinh^2 r}.
$$
Solving for $d\phi/ds$ in \thetag{6.2} and expressing everything in terms of $r$, $\rho$, and $\sin\a$, we obtain
$$
\der\phi s = \kg \frac{\cosh\rho - \sin\a\sinh\rho}{\cosh^2 r},
$$
and so
$$
\der\phi\theta = 
\kg \frac{\cosh\rho - \sin\a\sinh\rho}{\cosh^2r}\,\frac{\sinh^2\rho}{\sinh r}.
$$
This and the previous expression for $d\phi/d\theta$ in \thetag{6.1} are easily seen to be equal since $\sin\a = \cos\b = \tanh x/\tanh\rho$ by \thetag{A.3}.

In closing, we note that if $\phit$ is the central angle for the ``other end'' of the support line as in Figure~11, the symmetry in the relationship between $\phi$ and $\phit$ as related to $\omega$ in the proof of the corollary in \section4 appears as an asymmetry in the formulas for $d\phi/d\omega$, $d\phi/d\theta$, and $d\phi/ds$. For example,
$$
\der\phit\theta = 
\kg \frac{\cosh\rho + \sin\a\sinh\rho}{\cosh^2r}\,\frac{\sinh^2\rho}{\sinh r}
\qqand
\der\phit\omega = 1 + \frac{\sqrt{-k}}{c(x)}.
$$

\bigskip

\heading Appendix \endheading

In this appendix we recall formulas for $L(r)$, $A(r)$, and $c(r)$, and use them to express some standard trigonometric formulas in $\E2$, $\S2$, and $\H2$ in a unified manner. Let $M$ be $\S2$, $\E2$, or $\H2$ with curvature $k$. Let $L(r)$, $A(r)$, and $c(r)$ be the circumference, area, and curvature, respectively, of a circle of radius $r$ in $M$. It is convenient to define $\l(r) = L(r)/2\pi$ and $a(r) = A(r)/2\pi$. Familiar formulas for these are given in the following table.
$$
\spreadmatrixlines{10pt} 
\matrix
\format \c\qquad & &\qquad\c\qquad \\
          & k >0 & k =0 & k <0 \\
\ell(r) & \dfrac{\sin(\sqrt k r)}{\sqrt k } & r & \dfrac{\sinh(\sqrt{-k }r)}{\sqrt{-k }} \\
a(r) & \dfrac{1 - \cos(\sqrt k r)}k & r^2/2 & \dfrac{1 - \cosh(\sqrt{-k }r)}k \\
c(r) &\sqrt k\cot(\sqrt k r)&1/r&\sqrt{-k }\coth(\sqrt{-k }r)
\endmatrix \tag A.1
$$
It will also be convenient on occasion for $r$ to be negative, in which case $L(r)$ and $c(r)$ are taken to be negative as well. 

Note that $\l(r)c(r) = 1 - k a(r)$, which, when $k=0$, reduces to the familiar fact that the radius and curvature of a circle are reciprocals in Euclidean geometry. This quantity appears below and figures prominently in \section5. It is interesting to note that $\l(r)c(r)\Delta\theta = (1 - k a(r))\Delta\theta$ is an expression of the Gauss-Bonnet theorem for the circular sector of radius $r$ and angle $\Delta\theta$.

Consider a triangle in $M$ with sides $a$, $b$, $c$, and angles $\a$, $\b$, $\g$, with $a$ opposite $\a$, $b$ opposite $\b$, and $c$ opposite $\g$.
Using the functions defined above, the law of sines for all three geometries (see, {\sl e.g.}, \cite{R}) becomes
$$
\frac{\l(a)}{\sin\a} = \frac{\l(b)}{\sin\b} = \frac{\l(c)}{\sin\g}.
$$
According to Coxeter \cite{Co}, Bolyai noted that the law of sines could be expressed in this unified manner.
When $k\ne0$ the two laws of cosines in $\S2$ and $\H2$ become
$$
k\cos\g = \frac{\l(c)c(c) - \l(a)c(a)\,\l(b)c(b)}{\l(a)\l(b)} \qqand
\l(b)c(b) = \frac{\cos\b + \cos\a\cos\g}{\sin\a \sin\g}.
$$
These also hold when $k=0$, in which case the second one is equivalent to $\a + \b + \g = 0 \pmod{2\pi}$. Note that an appropriate limit of the first one yields the Euclidean law of cosines.

Now assume $\g = \pi/2$. The first law of cosines yields
$$
\l(c)c(c) = \l(a)c(a)\,\l(b)c(b), \tag A.2
$$
which is a version of the Pythagorean theorem when $k\ne0$. This, the law of sines, and the second law of cosines easily imply
$$
\spreadlines{.2truein}
\gathered
\sin\b = \frac{\l(b)}{\l(c)}, \quad \cos\b = \frac{c(c)}{c(a)}, \quad \tan\b = \frac1{\l(a)c(b)}, \qand \cos\b = \l(b)c(b) \sin\a.
\endgathered \tag A.3
$$
Note that these reduce to familiar Euclidean formulas when $k=0$.

In hyperbolic geometry we get special formulas for a triangle with one vertex at infinity. Note that $c(r)$ tends to $\sqrt{-k}$ as $r\to\infty$, which is the geodesic curvature of a horocycle. If the vertex at $\a$ recedes to a point at infinity while $a$ stays fixed, then $\b$ becomes the angle of parallelism for the length $a$, and we have
$$
\cos\b = \frac{\sqrt{-k}}{c(a)}, \qquad \sin\b = \frac1{\l(a)c(a)}, \qqand
\tan\b = \frac1{\sqrt{-k}\l(a)}.
\tag{A.4}
$$

\enddocument

\enddocument